\documentclass[11pt,reqno]{amsart}

\usepackage{url}
\usepackage[colorlinks=true, linkcolor=blue, citecolor=ForestGreen]{hyperref}
\usepackage{cite}

\usepackage{latexsym,amsmath,amssymb,amsfonts,amsthm}
\usepackage{mathrsfs}
\usepackage{mathtools}
\usepackage{dsfont}
\usepackage{bbm}
\usepackage{soul}     
\usepackage[utf8]{inputenc}
\usepackage{ifthen}

\usepackage[usenames, dvipsnames]{color}
\usepackage[usenames, dvipsnames, cmyk]{xcolor}
\usepackage{graphicx}
\usepackage[scriptsize,hang,raggedright]{subfigure}
\usepackage{multirow}
\usepackage{enumerate}
\usepackage{enumitem}
\usepackage{esint}

\usepackage[margin=1.2in]{geometry}
\usepackage{epsfig,epstopdf}

\mathtoolsset{showonlyrefs}

\numberwithin{equation}{section}
\definecolor{grey}{rgb}{.7,.7,.7}
\definecolor{refkey}{gray}{.45}
\definecolor{labelkey}{gray}{.45}

\newcommand{\xupref}[2]{\hspace{-0.3ex}\stackrel{\eqref{#1}}{#2}} % for (x)
\newcommand{\upupref}[3]{\hspace{-3ex}\stackrel{\eqref{#1},\eqref{#2}}{#3}}

\newtheorem{theorem}{Theorem}[section]
\newtheorem{proposition}[theorem]{Proposition}
\newtheorem{lemma}[theorem]{Lemma}

\newtheorem{corollary}[theorem]{Corollary}

\theoremstyle{remark}
\newtheorem{remark}[theorem]{Remark}
\theoremstyle{definition}

\newtheorem{example}[theorem]{Example}

\usepackage{tikz}
\usepackage{pgf,pgfplots}
\usetikzlibrary{calc}
\usetikzlibrary{decorations.markings}
\usetikzlibrary{decorations.pathmorphing}
\usetikzlibrary{decorations.shapes}
\usetikzlibrary{shapes,arrows,shapes.geometric,patterns,fadings}

\renewcommand{\a}{\mathbf{a}}
\renewcommand{\d}{\,\mathrm{d}}

\newcommand{\dx}{\d x}
\newcommand{\dy}{\d y}
\newcommand{\E}[1][]{
	\ifthenelse{\equal{#1}{ }}
	{\operatorname{\mathcal{E}_{\gamma}}}
	{\operatorname{\mathcal{E}_{\gamma}}(#1)}
}

\newcommand{\N}[1]{\operatorname{\mathcal{V}}(#1)}

\newcommand{\R}{\mathbb{R}}

\renewcommand{\div}{\operatorname{div}}
\newcommand{\spt}{\operatorname{spt}}

\newcommand{\Per}[1]{\operatorname{\mathcal{P}}(#1)}
\newcommand{\aPer}[1]{\operatorname{\mathcal{P}_{\a}}(#1)}

\newcommand{\mres}{\mathbin{\vrule height 1.6ex depth 0pt width
		0.13ex\vrule height 0.13ex depth 0pt width 1.3ex}}

\newcommand{\La}[1]{|#1|_{\a}}
\newcommand{\Ha}[1]{\mathcal{H}^{d-1}_{\a}(#1)}
\newcommand{\Ca}{C_{\a}}
\newcommand{\Ra}{R_{\a}}

\newcommand{\avgint}{\fint}

\newcommand{\e}{\varepsilon}

\renewcommand{\setminus}{\backslash}
\newcommand{\defeq}{\coloneqq}

\newcommand{\per}{\mathcal{P}}

\newcommand{\en}{\mathcal{E}_\gamma}

\newcommand{\ba}{\begin{array}}
\newcommand{\ea}{\end{array}}

\newcommand{\tld}[1]{\widetilde{#1}}
\newcommand{\bthm}{\begin{theorem}}
\newcommand{\ethm}{\end{theorem}}
\newcommand{\bcor}{\begin{corollary}}
\newcommand{\ecor}{\end{corollary}}
\newcommand{\bprop}{\begin{proposition}}
\newcommand{\eprop}{\end{proposition}}
\newcommand{\blemma}{\begin{lemma}}
\newcommand{\elemma}{\end{lemma}}
\newcommand{\bexmpl}{\begin{example}}
\newcommand{\eexmpl}{\end{example}}

\newcommand{\beqn}{\begin{equation}}
\newcommand{\eeqn}{\end{equation}}
\newcommand{\beqns}{\begin{equation*}}
\newcommand{\eeqns}{\end{equation*}}

\newcommand{\pr}{\prime}
\newcommand{\pt}{\partial}

\newcommand{\Rd}{\mathbb{R}^d}

\newcommand{\ol}{\overline}

\newcommand{\Sd}{\mathbb{S}^{d-1}}

\newcommand{\Hn}{\mathcal{H}^{d-1}}
\newcommand{\Hnm}{\mathcal{H}^{d-2}}

\newcommand{\V}{\mathcal{V}}

\renewcommand{\leq}{\leqslant}
\renewcommand{\geq}{\geqslant}

\definecolor{mygreen}{rgb}{0.1,0.75,0.2}

\newcommand{\varE}{\mathcal{E}}

\newcommand{\Om}{\Omega}

\newcommand{\calV}{\mathcal{V}}
\newcommand{\calF}{\mathcal{F}}

\newcounter{myenumi}
\DeclareMathOperator{\dive}{div}

\DeclareMathOperator*{\dist}{dist}
\DeclareMathOperator{\loc}{loc}

\DeclareMathOperator{\img}{Im}

\renewcommand{\le}{\leq}
\renewcommand{\ge}{\geq}

\title[A nonlocal isoperimetric problem with density perimeter]{A nonlocal isoperimetric problem with density perimeter}

\author{Stan Alama}
\address[Stan Alama]{Department of Mathematics and Statistics, McMaster University, Hamilton, ON, Canada}
\email{alama@mcmaster.ca}

\author{Lia Bronsard}
\address[Lia Bronsard]{Department of Mathematics and Statistics, McMaster University, Hamilton, ON, Canada}
\email{bronsard@mcmaster.ca}

\author{Ihsan Topaloglu}
\address[Ihsan Topaloglu]{Department of Mathematics and Applied Mathematics, Virginia Commonwealth University, Richmond, VA, USA}
\email{iatopaloglu@vcu.edu}

\author{Andres Zuniga}
\address[Andres Zuniga]{Department of Mathematics and Statistics, McMaster University, Hamilton, ON, Canada \& Instituto de Ciencias de la Ingenier\'{i}a, Universidad de O'Higgins (UOH), Rancagua, Chile  }
\email{andres.zuniga@uoh.cl}

\date{\today}                                        


\begin{document}

\begin{abstract}
We consider the minimization of an energy functional given by the sum of a density perimeter and a nonlocal interaction of Riesz type with exponent $\alpha$, under volume constraint, where the strength of the nonlocal interaction is controlled by a parameter $\gamma$. We show that for a wide class of density functions the energy admits a minimizer for any value of $\gamma$. Moreover these minimizers are bounded. For monomial densities of the form $|x|^p$ we prove that when $\gamma$ is sufficiently small the unique minimizer is given by the ball of fixed volume.  In contrast with the constant density case, here the $\gamma\to 0$ limit corresponds, under a suitable rescaling, to a small mass $m=|\Om|\to 0$ limit when $p<d-\alpha+1$, but to a large mass $m\to\infty$ for powers $p>d-\alpha+1$.
\end{abstract}

\maketitle
\tableofcontents

%%%%%%%%%%%%%%%%%%%%%%%%%%%%%%%%%%%%%%%%%%%%%%%%%%%%%%%%%%%%%%%%%
%%%%%%%%%%%%%%%%%%%%%%%%%%%%%%%%%%%%%%%%%%%%%%%%%%%%%%%%%%%%%%%%%

\section{Introduction}\label{sec:intro}

We consider the nonlocal isoperimetric problems
	\beqn \label{eq:nlip}
		e(\gamma) \defeq \inf \Big\{ \en(\Om) \colon |\Om|=1 \Big\}
	\eeqn
over sets of finite perimeter $\Om \subset \R^d$ with given volume, where $|\cdot|$ denotes the Lebesgue measure in $\R^d$, and the energy functional $\en$ is defined as
\beqn \label{eq:energy}
\en(\Om) \defeq \int_{\partial^*\Om} \a(x)\d\Hn + \gamma\int_\Om\int_\Om \frac{1}{|x-y|^\alpha}\d x\d y
\eeqn
for $\gamma>0$, $\alpha\in(0,d)$. Here $\pt^*$ denotes the reduced boundary of a set. The first term in the energy functional is the perimeter of $\Om$ with density $\a \colon \R^d\to[0,\infty)$, whereas the second term is a Riesz-type nonlocal interaction energy. 

The minimization problem \eqref{eq:nlip} is a variant of the classical liquid drop model introduced by Gamow in \cite{Ga1930}. Gamow's model is simply given by  \eqref{eq:nlip} with $\a\equiv 1$. The most important feature of this geometric variational problem is that the two terms present in the energy functional are in direct competition. For $\a \equiv 1$, the surface energy is minimized by a ball whereas the repulsive term does not admit a minimizer and prefers minimizing sequences with multiple vanishingly small components diverging infinitely apart in order to disperse the mass. The parameter of the problem, that is $\gamma$, sets a length scale between these competing forces and drives the competition between the short- and long-range interactions. This problem has generated considerable interest in the calculus of variations community (see e.g. \cite{AlBrChTo2017_3,BoCr14,ChPe2010,FFMMM,FKN16,FrLi2015,Ju2014,Julin2017,KnMu2014,KnMuNo2016,LuOtto2014,MurZal,RW2014} as well as \cite{ChMuTo2017} for a review) with several papers studying parameter regimes of existence and nonexistence of minimizers. Results of \cite{FKN16,KnMu2014,LuOtto2014}, for example, show that for large values of $\gamma$, the energy $\en$ with $\a\equiv 1$ does not admit a minimizer. There are also several studies characterizing the minimizing sequences \cite{AlBrChTo2017_3,BoCr14,KnMuNo2016} even when minimizers fail to exist. In particular, in \cite{AlBrChTo2017_3}, the authors use a ``regularization'' of the energy by adding an attractive external potential which guarantees the existence of minimizers for all values of $\gamma$.

Also very recently there has been studies on the extensions of the liquid drop model to the anisotropic setting where the surface energy is replaced by an anisotropic surface tension \cite{BoCrTo,ChNeuTo20,MiTo}. In these models the surface energy is given by $\int_{\pt^*\Om} \psi(\nu_\Om)\d\Hn$ for some convex, one-homogeneous function $\psi$ where $\nu_\Om$ denotes the outward unit normal to the reduced boundary $\pt^* \Om$. Such anisotropic extensions do not annihilate the translation invariance of the liquid drop model and a simple scaling argument heuristically justifies that for large $\gamma$ values minimizers still fail to exist. In contrast, the inclusion of a translation \emph{variant} density in the perimeter functional ``regularizes'' the liquid drop model in the sense that the problem admits a minimizer for all values of $\gamma$.

Isoperimetric problems defined via weighted perimeters
\beqn \label{eq:perimeter}
	\aPer{\Om}	\defeq \int_{\partial^*\Om}\a(x) \d\Hn(x)
\eeqn
have been studied for various choices of densities $\a$. Problems where the volume constraint is also weighted either by the density $\a$ or by some other function have especially attracted significant interest (see e.g. \cite{AlBrChMePo,barchiesibrancolinijulin,BrDePhilRu,BrChMe,CiGlPrRoSe20,CiPr2017,CaRoSe,DePhilFrPr,FuMaPr2011,MoPr,PrSa2018,PrSa,RoCaBaMo} and references therein). The main questions regarding these problems have been existence, boundedness and regularity of isoperimetric sets. These questions have been studied not only for specific densities (radial, monomial, Gauss-like) but also for rather general densities satisfying some boundedness and continuity conditions. To our knowledge, perturbations of density perimeter (either by long-range interactions or by external potentials) have not yet been considered in the literature.

\medskip

As pointed above, the inclusion of a confining density in the perimeter functional provides a different type of ``regularization'' of the problem \eqref{eq:nlip}. Our first main result establishes the compactness of any minimizing sequence with global convergence to a minimizer for any $\gamma\geq 0$ and for a wide class of densities, satisfying a simple coercivity condition: 
\medskip
		\begin{enumerate} [label=(\textbf{A\arabic*}),ref=(A\arabic*)]
				\addtolength{\itemsep}{8pt}
						\item \label{itm:A1} $\a \in C^0(\Rd)$, $\a(0)=0$, $\a(x)>0$ for all $x\neq 0$, and $\lim_{|x|\to \infty} \a(x) = +\infty$.
		\end{enumerate}
		\medskip

In order to state the existence result, we recall that sets $\Om_n\to \Om$ {\it globally} if $|\Om_n\triangle \Om|\to 0$, that is, their characteristic functions $\chi_{\Om_n}\to \chi_\Om$ in the $L^1(\R^d)$-norm.

\bthm[Existence of minimizers]\label{theorem:1}
	Let $\a$ be any density satisfying the assumption {\rm\ref{itm:A1}}, and fix any $\gamma\ge 0$.  Then any minimizing sequence $\{\Om_n\}_{n\in\mathbb{N}}$ for the  problem \eqref{eq:nlip} admits a subsequence which converges globally to a  minimizer  $\Om_\gamma\subset\R^d$ with $|\Om_\gamma|=1$ and  $ \chi_{\Om_\gamma}\in BV_{\loc}(\R^d\setminus\{0\})$.
\ethm

\begin{remark}\label{rem:vanishing}
	 The proof of the existence result in Theorem~\ref{theorem:1} can be extended to a somewhat broader class of densities $\a$ satisfying
	\[  
	\a\in C^0(\R^d),\quad \a^{-1}(\{0\})=\mathsf{Z}, \quad \a>0\,\text{ in }\,\R^d\setminus\mathsf{Z},\, \text{ and }\, \lim_{|x|\to\infty}\a(x)=+\infty,
	\]
	for any finite set $ \mathsf{Z} =\{z_1,\ldots,z_m\}\subset\R^d$ with $ m\in\mathbb{N}$.  We present the proof in the case $ \mathsf{Z}=\{0\} $. In this general situation, any minimizer $\Om_{\mathsf{Z}}\subset \R^d $ is such that $ \chi_{\Om_{\mathsf{Z}}}\in BV_{\loc}(\R^d\setminus\mathsf{Z}) $.
\end{remark}

We remark that the coercive nature of $\a$ at infinity ensures the existence of minimizers for \eqref{eq:nlip}, essentially because the splitting of mass off to infinity (the main reason for noncompactness in nonlocal isoperimetric problems) is rendered too costly.  However it does not ensure that minimizers need be connected sets.  Indeed, for large $\gamma$ the nonlocal interactions should become large enough to favor the fragmentation of sets, which will repel but be contained at finite distance. This behavior is also observed in nonlocal isoperimetric problems with a confining term \cite{AlBrChTo2017_3}.

Our existence result relies on a modified version of the relative isoperimetric inequality on annulli and requires only the minimal assumption {\rm\ref{itm:A1}} on the densities $\a$. On the other hand, proving boundedness of minimizers is rather technical and we prove it under one of the following additional structural assumptions:  
		\begin{enumerate} [label=(\textbf{A2\alph*}),ref=(A2\alph*)]
				\addtolength{\itemsep}{8pt}
						\item \label{itm:A2} $\a \in C^{0,1}_{\loc}(\Rd)$ and there exists constants $C_{\a} \geq 1$ and $R_{\a} \geq 1$ such that for every $R \geq R_{\a}$ it holds
																		\beqn\label{eq:growth}
																		0<\sup_{B_{2R}\setminus B_{R}}\a  \leq C_{\a} \inf_{B_{2R}\setminus B_{R}}\a.
																		\eeqn
				\item \label{itm:A3} $\a \in C^{0,1}_{\loc}(\Rd)$, $\a(x)=\a(|x|)$ and there exists $R_{\a} \geq 1$ such that $\a$ is non-decreasing for $|x| \geq R_{\a}$.
		\end{enumerate}	
\medskip

For densities satisfying either of the additional conditions {\rm\ref{itm:A2}} or {\rm\ref{itm:A3}}, using regularity results of quasi-minimizers, such as density bounds, we obtain the boundedness of minimizers of the problem \eqref{eq:nlip}.

\bthm[Boundedness of minimizers]\label{thm:bd}
 For any density $\a$ satisfying {\rm\ref{itm:A1}} and either {\rm\ref{itm:A2}} or {\rm\ref{itm:A3}}, and for any $\gamma\ge 0$, any minimizer $\Om_\gamma$ of~\eqref{eq:nlip} is essentially bounded.
\ethm

\medskip

\begin{remark}
We require the assumption \ref{itm:A1} in Theorem~\ref{thm:bd} only to obtain the existence of a minimizer to \eqref{eq:nlip}. If the existence of a minimizer could be obtained under some other conditions, either \ref{itm:A2} or \ref{itm:A3} would be sufficient to obtain the boundedness of minimizers.
\end{remark}

\begin{remark}[Almost polynomial densities]
Let $ \a\in C^{0,1}_{\loc}(\R^d) $ satisfy the following conditions:
	\begin{enumerate}[label=(\roman*)]
		\item $ \a(0)=0 $ and $ \a(x)>0 $ for $ x\neq 0 $.
		\item There exist $ p>0$, $ R_0>0 $ and $C_1,C_2>0 $ such that $C_1|x|^p \leq \a(x)\leq C_2|x|^p$ for all $|x|>R_0$.
	\end{enumerate}
Then $\a$ satisfies the assumptions \ref{itm:A1} and \ref{itm:A2}; hence, Theorems \ref{theorem:1} and \ref{thm:bd} hold for such densities.
\end{remark}

For homogeneous densities, which satisfy the condition $\a(tx) = t^p\, \a(x)$ for some $p>0$, 
a simple scaling argument shows that the minimization problem \eqref{eq:nlip} is equivalent to the problem
	\beqn\label{nlip1}
		 \inf \Big\{ \mathcal E_{1}(\Om) \colon |\Om|=m \Big\}
	\eeqn
with the correspondence
$$\gamma = m^{-(p+\alpha-d-1)/d}, \qquad p\neq p_* \defeq d-\alpha+1.  $$
It is interesting to observe that for homogeneous weights the large mass/small mass behavior of the minimization problem depends on the specific power $p$.  In particular, when $p>p_*$, the corresponding value of $\gamma$ varies inversely with mass.  Thus, with $p>p_*$
the nonlocal energy is dominated by the perimeter term $\per_{\a}$ for {\it large} mass $m$, and the nonlocal term dominates for {\it small} $m$, exactly the opposite of the behavior for constant $\a$.  For subcritical $p<p_*$ the opposite is true, and the energy is perimeter-dominated for small $m$.   
 At the critical value $p=p_*$ the two problems \eqref{eq:nlip} and \eqref{nlip1} are not equivalent, and \eqref{nlip1} is scale invariant:  the minimizers at any mass $m$ are all rescaled copies of the same set.
Since Theorem~\ref{theorem:1} guarantees the existence of a minimizer for all values of $m$ (or $\gamma$), an interesting question is the characterization of minimizers for a range of values of the parameters. We provide a partial answer to this question in the next theorem.

\bthm[Global minimizers in the small $\gamma$ regime] \label{thm:global}
Let $\a(x)=|x|^p$ with $p>0$. For $\gamma$ sufficiently small the ball $B\subset\Rd$ of volume one, centered at the origin is the unique minimizer of $e(\gamma)$.
\ethm 

For $p<p_*$, the $\gamma\to 0$ limit is equivalent to the small $m$ regime in \eqref{nlip1}, and the optimality of the spherical ball for small mass is well-known for the unweighted $\a\equiv 1$ case (see \cite{BoCr14,Ju2014,Julin2017,KnMu2013,KnMu2014}). With $p>p_*$, the situation is reversed and Theorem~\ref{thm:global} shows that the ball minimizes for {\it all sufficiently large} $m$.  Such $\a$ are very coercive at infinity, and the situation is similar to the case studied by G\'{e}n\'{e}raux and Oudet in \cite{GenOud}, where the problem \eqref{nlip1} with constant density is supplemented with a confining term. In that case they also prove the minimality of the ball for very large $m$.

The proof of Theorem~\ref{thm:global} relies on a penalization technique, similar to those used in several geometric variational problems involving the perimeter functional (see \cite{AcFuMo13,BoCr14,CicaleseLeonardi,FFMMM,Neum16}). Utilizing results from the regularity theory for density perimeters  \cite{CiPr2017,PrSa} we reduce the minimizers of the nonlocal problem to nearly spherical or isoperimetric sets in the small $\gamma$ regime. The novelty here is, though, that we cannot directly apply the results from the literature due to the degeneracy of the density $\a$ at the origin and the possibility of small nonsmooth components of $\pt^*\Om_\gamma$ near the origin. Once we reduce the problem to nearly spherical sets we use a Fuglede-type argument (see \cite{Fuglede}) to control the isoperimetric and nonlocal deficits between minimizers and the ball and show that for small $\gamma>0$ these quantities have to be identically zero.

The regime of large $\gamma$ is also very interesting, but its analysis requires a very different approach. While the existence of minimizers is guaranteed by Theorem~\ref{theorem:1}, in this regime, the nonlocal term is dominant and prefers the minimizer to break into smaller pieces distributed in a compact set whose size is determined by the confining term $\a$. Hence, the characterization of minimizers (i.e., the shape of the disconnected components as well as their locations) depends on the delicate balance between the preferred shapes dictated by the density perimeter and the inter-component interactions. A similar phenomenon is also observed in models of copolymer/homopolymer blends. While existence of minimizers is obtained for all values of $m$ in \cite{BoKn2016}, only in the small $m$ regime the minimizers are uniquely characterized leaving the question of the precise morphology of minimizing configurations for large $m$ open.

\addtocontents{toc}{\protect\setcounter{tocdepth}{1}}
\subsection*{Structure of the paper}
The paper is organized as follows. In Section~\ref{sec:exist} and Section~\ref{sec:bounded} we prove Theorem~\ref{theorem:1} and Theorem~\ref{thm:bd}, respectively. Section~\ref{sec:globalmin} is devoted to the proof of Theorem~\ref{thm:global}. 

\subsection*{Notation}
Throughout the paper $\omega_d$ denotes the volume of the unit ball $ B_1(0)$ in $ \R^d $ and we write $ B_r \defeq B_r(0)$ to denote the ball of radius $r$ centered at zero. Constants, denoted by $C$, can change from line to line (unless otherwise noted). We will denote the $\a$-volume measure and the $\a$-surface area measure, respectively, by
\[  
	\La{\Om}\defeq \int_{\Om}\a(y)\dy, \quad \text{ and } \quad \Ha{\Om} \defeq \int_{\Om}\a(y)\,\d\Hn(y).
\]
The relative weighted perimeter of $E$ in $F$ will be denoted by either $\aPer{E,F}$ or $ \int_F \a(x)|\nabla \chi_E|$,
where $\chi_E$ is the characteristic function of the set $E$ and $|\nabla \chi_E|$ is the total variation of $\chi_E$. The weighted perimeter of $E$ in $F$ (or the $\a$-perimeter) is defined as
	\[
		\aPer{E,F} :=\sup \left\{ \int_{E} \dive\big(\a(x) X(x)\big) \dx \colon X\in C_c^\infty(F;\R^d), \ \|X\|_{L^\infty} \leq 1 \right\}.
	\]
In particular, $\aPer{E,F}=\Ha{\pt^* E \cap F}$. We will say that $\chi_E \in BV_{\a}(F)$ if $\aPer{E,F} < +\infty$.
Perimeters of sets in the whole space (i.e., when $F=\Rd$) are denoted by $\mathcal{P}_{\a}(E)$ or $\int_{\Rd} \a(x)|\nabla \chi_E|$. We will denote the Euclidean perimeter (when $\a \equiv 1$) by simply $\mathcal{P}$. Finally, we will denote the nonlocal term by $\mathcal{V}$, i.e.,
\[  
	\N{\Om} \defeq \int_\Om \int_\Om \frac{\dx\dy}{|x-y|^{\alpha}},
\]
for any $\alpha \in (0,d)$. 
\addtocontents{toc}{\protect\setcounter{tocdepth}{2}}

%%%%%%%%%%%%%%%%%%%%%%%%%%%%%%%%%%%%%%%%%%%%%%%%%%%%%%%%%%%%%%%%%
%%%%%%%%%%%%%%%%%%%%%%%%%%%%%%%%%%%%%%%%%%%%%%%%%%%%%%%%%%%%%%%%%

\section{Existence of minimizers}\label{sec:exist}

In this section we present the proof of Theorem~\ref{theorem:1}. It relies on the following modified version of the relative isoperimetric inequality on annulli.

\blemma\label{lem:isop:annuli}
	Let $ A_{r,R}=\big\{x\in\R^d: r\leq |x|< R\big\}=B_R\setminus B_r $ denote the half-open annulus of inner radius $ r $ and outer radius $ R $. Then, there exists a dimensional constant $ c_d>0 $, so that 
	\begin{equation}\label{eq:lem:relativeisop}
	\min\Big\{|\Om\cap A_{1,2}|^{(d-1)/d},\, |A_{1,2}\setminus \Om|^{(d-1)/d}\Big\}
	\leq c_d\,\mathcal{P}(\Om,A_{1,2})
	\end{equation}
	for every set of finite perimeter $ \Om \subset\R^d $.
\elemma
	
\begin{remark}\label{rem:rel-isop}
	Inequality~\eqref{eq:lem:relativeisop} is still valid with the same constant $ c_d $ over any annulus $ A_{2^j,2^{j+1}} $ with $ j\geq 1 $, as the inequality is invariant under scalings.
\end{remark}
The relative isoperimetric inequality is typically stated for balls in $\R^d$ (see e.g.~\cite[Thm 5.4.3]{ziemerWeakly},) but in fact the same proof verifies that it holds in any domain for which one can prove the validity of the Poincar\'e inequality,
\[  
\left(\int_{A_{1,2}}|u-\bar{u}_{r,R}|^{\frac{d}{d-1}}\,\dx\right)^{\frac{d-1}{d}}\leq
C_{1,2}\int_{A_{1,2}}|\nabla u|\dx ,
\]
where $ \bar{u}_{r,R} \defeq \avgint_{A_{r,R}}u\dx $ is the average on annuli. The latter can be found in~\cite{adamsSobolev}.

\bigskip

We now turn to the proof of the existence of minimizers for $e(\gamma)$.  For translation-invariant nonlocal isoperimetric problems existence is a delicate issue, as minimizing sequences can split, with pieces diverging to infinity.  The increasing weight $\a(x)$ raises the cost of splitting, an effect which is quantified in our proof via the relative isoperimetric inequality, Lemma~\ref{lem:isop:annuli}.

\begin{proof}[Proof of Theorem \ref{theorem:1}] Let $ \{\Omega_n\}\subset \R^d$ be a minimizing sequence of~\eqref{eq:nlip}:
\[  
	\en(\Omega_n)\to e(\gamma)\;\text{ as }\;n\to\infty,\;\;\text{ and }\;\; |\Omega_n|=1\;\;\text{for any }n\geq 1.
\]
We first show that $ \{\Omega_n\} $ is uniformly bounded in $ BV(\R^d\setminus \overline{B}_{\e}) $ for all $ \e>0 $ sufficiently small. To see this, note that $ \lim_{|x|\to\infty}\a(x)=+\infty $ implies that there exists $ R_1>0 $ such that $ \a\geq 1 $ in $ \R^d\setminus B_{R_1} $. Define $ \delta_{\e} \defeq \min\{\a(x) \colon \e\leq |x|\leq R_1\} $ for $ 0<\e<1 $ small enough so that $ \delta_{\e}\leq 1 $ (as $ \lim_{|x|\to 0 }\a(x)=0 $). Then, for any such choice of $ \e $, $ \a(x)\geq \delta_{\e}>0 $ for all $ x\in\R^d\setminus \ol{B}_{\e} $, and we deduce
\[  
\delta_{\e}\,\Per{\Omega_n,\R^d\setminus\overline{B}_{\e}}<\aPer{\Omega_n,\R^d\setminus \overline{B}_{\e}}<\en(\Omega_n)=e(\gamma)+o_n(1),
\]
which confirms uniform boundedness in $ BV(\R^d\setminus \overline{B}_{\e}) $.

We next show local convergence of $\{\Omega_n\}$ to a limiting set $\Om_\gamma$.
Invoking compactness results of sets with uniformly bounded perimeter, there exists $ \Omega^{\e}\subset\R^d $ so that $\Omega_{n_{\ell}}\to \Omega^{\e}$ locally in $\R^d\setminus \overline{B}_{\e} $, as $\ell\to \infty$.
Running a diagonalization argument over a sequence $ \e_k=1/k\to 0^+ $, there exists a subsequence ${n_\ell}\to\infty $ such that for all $k\geq 1 $,  $ \Omega_{n_\ell} \to \Omega^{1/k} $ locally in $ \R^d\setminus \overline{B}_{1/k} $, as $ \ell\to\infty $. In particular, $ 
\Omega^{1/(k+j)}\setminus \overline{B}_{1/k}=\Omega^{1/k}\setminus \overline{B}_{1/k} $ for any $ j\geq 1 $. Defining the limit set as $ \Om_\gamma \defeq \bigcup^{\infty}_{k=1} (\Omega^{1/k}\setminus\overline{B}_{1/k})$, we claim that up to subsequence, 
\beqn
	\begin{gathered}\label{eq:pf:thm1:1}
 \Omega_n\to \Omega_\gamma \;\text{ locally in }\R^d\;\text{ as }\; n\to \infty,\text{ and } \\
  \chi_{\Omega_n}\to \chi_{\Omega_\gamma} \;\text{ pointwise a.e. in }\R^d\;\text{ as }\; n\to \infty.
  	\end{gathered}
\eeqn
Assuming the claim, we may conclude that $ \Per{\Om_\gamma,\R^d\setminus\overline{B}_{\e}}\leq \liminf\limits_{n\to\infty}\Per{\Omega_n,\R^d\setminus\overline{B}_{\e}}$ for $ 0<\e \ll 1 $, which shows $ \chi_{\Om_\gamma}\in BV_{\loc}(\R^d\setminus \{0\}) $.

\medskip

To verify~\eqref{eq:pf:thm1:1}, let $ K\subset\R^d $ be a compact set and for $ \e>0 $, fix $ k \ge 1 $ so that $ |B_{1/k}|\leq \e $. Then
	\begin{align*}
			 |(\Omega_{n}\triangle\Om_\gamma)\cap K| &= |(\Omega_{n}\triangle\Omega^{1/k})\cap (K\setminus B_{1/k})|+ |(\Omega_{n}\triangle\Om_\gamma)\cap K\cap B_{1/k}| \\
																		 &\leq o_n(1)+\e.
	\end{align*} 
Since $ \e>0 $ is arbitrary small, we conclude~\eqref{eq:pf:thm1:1}. 
The aforementioned convergence along with the lower semicontinuity of the $ \a $-perimeter functional, together with Fatou's Lemma shows that
\[  
	\en(\Om_\gamma)=\aPer{\Om_\gamma}+\gamma\N{\Om_\gamma}\leq \liminf _{n\to\infty}(\aPer{\Omega_n}+\gamma \N{\Omega_n})=e(\gamma).
\]
We are only left to show that $ \Om_\gamma $ is admissible in~\eqref{eq:nlip}, from which it will follow that $ \en(\Om_\gamma)\geq e(\gamma) $; thus obtaining the existence of a minimizer. 

\medskip

We observe that $ |\Om_\gamma|\leq 1 $ in view of Fatou's Lemma, once again. We claim that in fact $ |\Om_\gamma|=1  $. Suppose, on the contrary, that $ |\Om_\gamma|<\beta $ for some $ \beta\in(0,1) $. 

The local convergence~\eqref{eq:pf:thm1:1} shows that for all $ R>0 $, $ |\Omega_n\cap B_R|=|\Om_\gamma\cap B_R|+o_n(1)<\beta $ for all but finitely many $ n $. 
Thus, the sets $\Om_n$ have very thick tails, which will introduce huge energy cost via the relative isoperimetric inequality.
By running a diagonalization argument over $ \{R_k=2^k\} $, there exists an increasing subsequence $ n_k\to+\infty $ such that for all $k\geq 1$,
\begin{equation}\label{eq:pf:thm1:2} 
	\inf_{n\geq n_k}|\Omega_n\setminus B_{2^k}|> 1-\beta.
\end{equation}
On the other hand, as $ \a(x)\to+\infty$ as $ |x|\to\infty $, for $ M>1 $ arbitrarily large, for all $ n\geq 1 $, and every $ k\geq k_M$ sufficiently large,
\begin{equation}\label{eq:pf:thm1:3}  
	\en(\Omega_n)\geq \aPer{\Omega_n,\R^d\setminus B_{2^k}} \geq M \Per{\Omega_n,\R^d\setminus B_{2^k}}.
\end{equation}

Intuitively, when $j$ is large we expect $|\Om_n\cap A_{2^{j},2^{j+1}}|$ to be much smaller than
its complement $|A_{2^j,2^{j+1}}\setminus \Omega_n|$ in the annulus.  Indeed,
we claim that there exists $j_0\in\mathbb{N}$ such that for all $j\geq j_0$ and for all but finitely many $ n $ we have
\begin{equation}\label{eq:pf:thm1:4}   
|\Omega_n\cap A_{2^{j},2^{j+1}}|<(2^{j+1})^{-d}|A_{2^j,2^{j+1}}\setminus \Omega_n|.
\end{equation}
For otherwise, there would exist increasing sequences 
$ j_{\ell}\to+\infty $ and $ n_{\ell}\to+\infty $ for which $ |\Omega_{n_{\ell}}\cap A_{2^{j_{\ell}},2^{j_{\ell}+1}}|\geq (2^{j_{\ell}+1})^{-d}|A_{2^{j_{\ell}},2^{j_{\ell}+1}}\setminus \Omega_{n_{\ell}}| $, for every $ \ell\geq 1 $.  This would imply that
\begin{equation*}  |\Omega_n\cap A_{2^{j},2^{j+1}}| 
             \ge {\frac{1}{1+2^{(j+1)d}}} |A_{2^{j_{\ell}},2^{j_{\ell}+1}}|, 
\end{equation*}
and hence,
\begin{align*}
	1=|\Omega_{n_{\ell}}|
	&\geq\sum\limits^{\infty}_{\ell=1}|\Omega_{n_{\ell}}\cap A_{2^{j_{\ell}},2^{j_{\ell}+1}}| \\
	&\geq\sum\limits^{\infty}_{\ell=1}\frac{1}{1+(2^{j_{\ell}+1})^d}|A_{2^{j_{\ell}},2^{j_{\ell}+1}}| >\sum\limits^{\infty}_{\ell=1}\frac{1}{2}\biggl(1-\frac{1}{2^{d}}\biggr)\omega_d=+\infty,
\end{align*}
establishing the claim.

We now fix any $ k\geq \max\{k_M,j_0\} $ and $ n \gg 1 $ sufficiently large to obtain the validity of~\eqref{eq:pf:thm1:3} and of~\eqref{eq:pf:thm1:4} for all $ j\geq k $. Utilizing the relative isoperimetric inequality on $ \{A_{2^j,2^{j+1}}: j\geq k\} $ (see Remark~\ref{rem:rel-isop}), we get the lower bound 
\begin{align*}
	c_d\, \Per{\Omega_n,\R^d\setminus B_{2^{k}}} 
	&=\sum^{\infty}_{j=k}c_d\, \Per{\Omega_n,A_{2^j,2^{j+1}}} \\
	&\geq \sum^{\infty}_{j=k}|\Omega_n\cap A_{2^j,2^{j+1}}|^{\frac{d-1}{d}}\\
	&\geq\biggl(\sum^{\infty}_{j=k}|\Omega_n\cap A_{2^j,2^{j+1}}|\biggr)^{\frac{d-1}{d}}=|\Omega_n\cap (\R^d\setminus B_{2^k})|^{\frac{d-1}{d}}
\end{align*}
Increasing the value of $ n\geq n_k $ if necessary, it follows from~\eqref{eq:pf:thm1:2} and \eqref{eq:pf:thm1:3} that
\[  
	+\infty>e(\gamma)+o_n(1) \geq M \Per{\Omega_n,\R^d\setminus B_{2^{k}}} \geq \frac{M}{c_d}(1-\beta)^{\frac{d-1}{d}},
\]
with $ M>1 $ arbitrarily large. Thus, we reach a contradiction. Hence, $ |\Om_\gamma|=1 $, and we have proven that $\Om_\gamma$ attains the minimum in the nonlocal isoperimetric problem \eqref{eq:nlip}.

Finally, by the identity $ |\Omega_n\triangle\Om_\gamma|= 2|\Om_\gamma\setminus \Omega_n|+|\Omega_n|-|\Om_\gamma|$ together with $ |\Omega_n|= 1=|\Om_\gamma|$ and~\eqref{eq:pf:thm1:1}, we deduce the global convergence of the subsequence $\{\Omega_n\}$, as (by local convergence,)
\[  
	|\Om_\gamma\setminus \Omega_n|\leq |(\Om_\gamma\setminus\Omega_n)\cap B_r|+|\Om_\gamma\setminus B_r|\leq |(\Om_\gamma\setminus \Omega_n)\cap B_r|+o_r(1).
\]
Thus, every minimizing sequence $\Om_n$ for $e(\gamma)$ contains a subsequence which converges globally (in $L^1(\R^d)$) to a minimizer of $e(\gamma)$.
\end{proof}

%%%%%%%%%%%%%%%%%%%%%%%%%%%%%%%%%%%%%%%%%%%%%%%%%%%%%%%%%%%%%%%%%
%%%%%%%%%%%%%%%%%%%%%%%%%%%%%%%%%%%%%%%%%%%%%%%%%%%%%%%%%%%%%%%%%

\section{Boundedness of minimizers}\label{sec:bounded}

In this section we prove Theorem~\ref{thm:bd}. Since the assumptions \ref{itm:A2} and \ref{itm:A3} characterize different types of densities, the proof of the theorem requires two different approaches. For densities which are polynomial-like and have bounded oscillations (i.e., densities satisfying \ref{itm:A2}) we make use of a series of technical lemmas establishing uniform density bounds for quasi-minimizers of the weighted perimeter functional where the density is measured with respect to weighted volumes. For radial and monotone densities (i.e., densities satisfying \ref{itm:A3}), on the other hand, we utilize a regularity result, called $\e-\e^{(d-1)/d}$ property, which basically says that a set of finite perimeter can be locally modified where one increases its volume by $\e$, while the perimeter increases at most by a constant multiple of $\e^{(d-1)/d}$. We present the proof in two subsections.

\addtocontents{toc}{\protect\setcounter{tocdepth}{1}}
\subsection{Densities with bounded oscillations} 

We start with densities satisfying the assumptions \ref{itm:A1} and \ref{itm:A2}. First, we prove that any minimizer $\Om_\gamma$ of \eqref{eq:nlip} has finite $\a$-volume.

\blemma
For any density $\a$ satisfying {\rm\ref{itm:A1}} and {\rm\ref{itm:A2}} any minimizer $\Om_\gamma$ of~\eqref{eq:nlip} has finite $\a$-volume, i.e., $|\Om_\gamma|_{\a}<+\infty$.
\elemma

\begin{proof}
  By passing to the limit $ \Omega_n\to\Om_\gamma $ in~\eqref{eq:pf:thm1:4}, 
\begin{equation}\label{eq:ann0}
|\Om_\gamma\cap A_{2^{j},2^{j+1}}|\leq (2^{j+1})^{-d}|A_{2^j,2^{j+1}}\setminus \Om_\gamma|,
\end{equation}
for all $j\ge j_0$, for some $j_0$.  Fix $ j_{\a}\in\mathbb{N} $ with $j_{\a}\ge j_0$, such that $ R_{\a}<2^{j_{\a}}$ where $R_{\a}$ is given as in \eqref{eq:growth}.
Using ~\eqref{eq:growth}, \eqref{eq:ann0}, and the relative isoperimetric inequality~(see Remark~\ref{rem:rel-isop})  we have:
\begin{align*}
	\int_{\Om_\gamma\setminus B_{2^{j_{\a}}}}\a(x)\dx
	&=\sum^{\infty}_{j=j_{\a}}\int_{\Om_\gamma\cap A_{2^j,2^{j+1}}}\a(x)\dx\\
	&\leq\sum^{\infty}_{j=j_{\a}}\big(\sup_{A_{2^j,2^{j+1}}}\a \big)|\Om_\gamma\cap A_{2^j,2^{j+1}}|\\
	&\leq\sum^{\infty}_{j=j_{\a}}\Ca\big(\inf_{A_{2^j,2^{j+1}}}\a\big)|\Om_\gamma\cap A_{2^j,2^{j+1}}|^{\frac{d-1}{d}}|\Om_\gamma\setminus B_{2^{j}}|^{\frac{1}{d}}\\
	&\leq \Ca\,|\Om_\gamma\setminus B_{2^{j_{\a}}}|^{\frac{1}{d}}\sum^{\infty}_{j=j_{\a}}\big( \inf_{A_{2^j,2^{j+1}}}\a \big) c_d\int_{A_{2^j,2^{j+1}}}|\nabla\chi_{\Om_\gamma}|\displaybreak[1]\\
	&\leq c_d\,\Ca\,|\Om_\gamma\setminus B_{2^{j_{\a}}}|^{\frac{1}{d}}\sum^{\infty}_{j=j_{\a}}\int_{A_{2^j,2^{j+1}}}\a(x)|\nabla\chi_{\Om_\gamma}|.
\end{align*}
Hence, as $ |\Om_\gamma|\leq 1 $, we conclude that
\begin{equation}\label{eq:pf:thm1:6}
	\int_{\Om_\gamma\setminus B_{2^{j_{\a}}}}\a(x)\dx\leq  c_d\,\Ca\;\int_{\R^d\setminus B_{2^{j_{\a}}}}\a(x)|\nabla\chi_{\Om_\gamma}|.
\end{equation}  
Since $\Om_\gamma$ has finite $\a$-perimeter, we obtain that $|\Om_\gamma|_{\a}<+\infty$.
\end{proof}

At the heart of the proof of Theorem \ref{thm:bd} lies the regularity of quasi-minimal sets with a volume constraint. In order to establish this we will largely follow the argument carried out by Rigot in~\cite[Chapter~2]{rigot2000ensembles}, where the author studies the case of standard perimeter functional $ \a(x)\equiv 1 $.  As in \cite[Chapter~1]{rigot2000ensembles} and~\cite[Chapter~21]{maggiBOOK}, given a function $ g:(0,+\infty)\to (0,+\infty)$ with $ g(x)=o(x^{(d-1)/d}) $ for $ x $ close to 0, we will say that $ \Om_\gamma $ is a \emph{volume constrained quasi-minimal set} for $ \a $-perimeter if 
\[  
	\aPer{\Om_\gamma} \leq \aPer{F}+g(|F\triangle \Om_\gamma|)
\]
for any $ F\subset\R^d $ with $ \chi_F\in BV_{\a}(\R^d)$, $ |F|=1 $ and $ F\triangle \Om_\gamma\subset\!\subset\R^d $.

Minimizers of isoperimetric problems with a Riesz-type nonlocal term are also volume constrained quasi-minimizers for the perimeter functional with the choice $ g_{\gamma}(x)\simeq \gamma \,x $. Indeed, as argued in~\cite[Proposition 2.1]{KnMu2014}, we define the potential of a Borel set $ E\subset\R^d $ by 
\[
	v_E(x) \defeq \int_{E}\frac{1}{|x-y|^{\alpha}}\dy
\] 
and consider a minimizer $ \Om_\gamma $ of~\eqref{eq:nlip} together with a set $ F $ with prescribed mass and $ F\triangle \Om_\gamma\subset B_r(0) $ for some $ r>0 $. A simple argument shows that the interaction energy $\calV$ is Lipschitz with respect to symmetric difference:
\beqn
	\begin{aligned} \label{eq:V_Lip}
		\N{F}-\N{E} &=\int_{F}\int_{F}\frac{1}{|x-y|^{\alpha}}\dx\dy-\int_{E}\int_{E}\frac{1}{|x-y|^{\alpha}}\dx\dy\\
										&\leq \int_{E\triangle F}(v_{E}+v_{F})\dx\\
										&\leq C |E\triangle F|
	\end{aligned}
\eeqn
with $C=2\int_{B_1}\frac{\dy}{|y|^{\alpha}}+2$.

Hence, any minimizer $ \Om_\gamma $ of $ \varE_\gamma$ must satisfy
\begin{equation}\label{eq:quasi-minimality}
	\aPer{\Om_\gamma}	\leq \aPer{F} + C\gamma\, |F\triangle \Om_\gamma|,
\end{equation} 
for any suitable competitor $ F $ as above.

We will now present some technical lemmas essentially studied in~\cite[Chapter~2]{rigot2000ensembles} and adapt these results to the case of weighted perimeters. The next lemma, proven for $\a=1$ by Giusti (see \cite[Lemma 2.1]{G}), shows that any set of positive perimeter can be approximated in $L^1$ by another set without substantially increasing the weighted perimeter.  We denote by $ [\, \cdot\,]_{1,D} $  the Lipschitz seminorm in $ D $.

\blemma\label{lem:Rigot:2.1.1}
	Let $ D\subset\R^d $ be a bounded domain, $ L\subset\R^d $ with $ \chi_L\in BV_{\a}(\R^d)$ for a density $ \a\in C^{0,1}(D) $, $ \inf_{D}\a>0 $, and such that
	\[  
		\int_D\a(x)|\nabla\chi_L|>0.
	\]
	Then there exist $ \e>0 $ and $ Q_{\a}>0 $, depending on $ L\cap D $ and $ D $ only, such that,
	$ Q_{\a}\lesssim \big(\sup_{D}\a\big) \big(1+[\a]_{1,D}\big) $, and
	 for every $ v\in (-\e,\e),$ there exists $F\subset\R^d $ with $ F=L $ in a neighborhood of $ \R^d\setminus D $ satisfying
	\begin{gather*}
	 |F|=|L|+v, \\
	 \int_D\a(x)|\nabla\chi_F|\leq \int_D\a(x)|\nabla\chi_L|+Q_{\a}|v|,\\
	 \int_D|\chi_F-\chi_L| \leq Q_{\a}|v|.
	\end{gather*}
\elemma

\begin{proof}
By definition of $ \a $-perimeter, there exists $ w\in C^1_c(D;\R^d) $, $ |w(x)|\leq \a(x) $ a.e. in $ D $, such that
\begin{equation} \label{eq:lem1Rigot:0}
	\int_{D}\chi_L\div w\dx\geq \frac{1}{2}\int_{D}\a(x)|\nabla\chi_L|>0.
\end{equation}
Note that, since $ \a>0 $ in $ D $,  $ BV_{\a}(D)\subset BV(D) $. For $ t\in(0,1) $ we put $ \eta_t=x+tw(x) $ and $ K \defeq \spt w\subset\!\subset D $. Then $ \eta_t\equiv I $ in $ \R^d\setminus K $, and for $ |t| $ small enough, $ \eta_t \colon D\to D $ is a diffeomorphism. Letting $ L_t \defeq \eta_t(L)$, we claim:
\beqn\label{eq:lem1Rigot:1}
\begin{gathered}
	|L_t|=\int_L|\det D\eta_t|\dx\\
	\int_{D}\a(x)|\nabla\chi_{L_t}|
	\leq \int_{D}\a(x)f_t(x)|\nabla\chi_L|+|t| \big(\sup_K\a \big)[\a]_{1,K}\int_{K}f_t(x)|\nabla\chi_{L}|
\end{gathered}
\eeqn
where $ f_t(x) \defeq |\det D\eta_t(x)||(D\eta_t)^{-1}(x)| $.  Indeed, the first equality is clear, and the inequality below it is obtained by noting,
\begin{align*}
	\int_{D}\a(x)|\nabla\chi_{L_t}|
	&\leq \int_{D}(\a\circ\eta_t)(x)\,f_t(x)|\nabla\chi_L| \\
	&\leq \int_{D}\a(x)\,f_t(x)|\nabla\chi_L|+\int_{K}f_t(x)\big|\a\circ\eta_t(x)-\a(x)\big|\,|\nabla\chi_L|
\end{align*}
and estimating, for $ x\in K $,
\[ 
	|\a\circ\eta_t(x)-\a(x)| = |\a(x+tw(x))-\a(x)| \leq [\a]_{1,K}\,|tw(x)|\leq |t|\, [\a]_{1,K}\big(\sup_{K}\a\big).
\]
Also, $ \det D\eta_t=1+t\div w+t^2A(x,t) $ and $ (D\eta_t)^{-1}=I-tH(x,t) $, with $ |A| $ and $ |H| $ bounded uniformly by a constant, which depends exclusively on $ L\cap D $ and $ D $. For $ |t| $ small enough, $ f_t(x)\leq 1+t(\div w+|H|)+O(t^2) $, and so~\eqref{eq:lem1Rigot:1} shows
\begin{equation}\label{eq:lem1Rigot:2}
	|L_t|=|L|+t\int_{D}\chi_L\div w\dx+t^2\int_{D}\chi_LA(x,t)\dx,
\end{equation}
and
\beqn \label{eq:lem1Rigot:3}
\begin{aligned}
	\int_{D}\a(x)|\nabla\chi_{L_t}|
	&\leq \big(1+|t|\cdot \|\div w+|H|\|_{L^{\infty}(D)}\big)\int_{D}\a(x)|\nabla\chi_L| \\
	&\qquad\qquad\qquad\qquad + |t| [\a]_{1,D} \big(\sup_{D}\a \big) \int_{D}|\nabla\chi_L|+O(t^2)\\
	&\leq \int_{D}\a(x)|\nabla\chi_L| +|t| \big(\sup_{D}\a \big) (C+\,[\a]_{1,D})\int_{D}|\nabla\chi_L|+O(t^2).
\end{aligned}  
\eeqn
In view of~\eqref{eq:lem1Rigot:0}, there exists $ \e^\pr>0$ sufficiently small so that for any choice of $ v\in(-\e',\e') $ the relation $ t\,\int_{D}\chi_{L}\div w\dx+t^2\int_{D}\chi_{L}A(x,t)=|v| $ in~\eqref{eq:lem1Rigot:2} holds true for some $ t_v $, and moreover $ |t_v|\leq C'|v| $, with $ C' $ depending on $ L\cap D $ and $ D $ only. We take $ F \defeq L_{t_v} $ and observe that $ F $ satisfies the first two statements of the lemma, in light of~\eqref{eq:lem1Rigot:2}-\eqref{eq:lem1Rigot:3}, for the value $ Q_{\a}=2(\sup_{D}\a)(C+[\a]_{1,D})\int_{D}|\nabla\chi_{L}| $, by decreasing the value of $ \e' $ if necessary, {\color{magenta} in order that $C'\e'\le 2$.}

 To verify the final statement, for $ g\in C^1(D) $ let $ g_t \defeq g\circ\eta^{-1}_t $, so $  
g_t-g = g_t-g_t\circ \eta_t$. Then
\begin{align}
	\int_{D}|g_t-g|\dx
	&=-\int_{D}\int^1_0tw(x)\cdot \nabla g_t(x+tsw(x))\d s\dx\notag\\
	&\leq|t|\int^1_0\int_{D}\a(x)|\nabla g_t\circ\eta_{ts}|\dx\d s\notag\\
	&\leq |t|\int_{D}\a(x)|\nabla g|\dx+|t|^2\,[\a]_{1,D}\int_{D}|\nabla g|+O(|t|^3),\label{eq:lem1Rigot:4}
\end{align}
where the last inequality will be derived below. Observe the third bound in the statement of Lemma~\ref{lem:Rigot:2.1.1} holds for $ Q_{\a}=2(\sup_D\a)(1+[\a]_{1,D})\int_{D}|\nabla\chi_{L}| $, upon decreasing the value of $ \e' $ if necessary. An approximation argument justifies estimate~\eqref{eq:lem1Rigot:4} for $ g\in BV_{\a}(\R^d) $ and so in particular for $ g=\chi_L $ and $ g_t=\chi_{L}\circ\eta^{-1}_t=\chi_{L_t} $. First, note that
\begin{align*}
	\int_{D}\a(x)|\nabla g_t\circ\eta_{ts}|\dx
	&=\int_{D}(\a\circ\eta^{-1}_{ts})(x)|\nabla g_t|\,|\det D(\eta^{-1}_{ts})|\dx\\
	&\leq (1+|t|\,\|H\|_{\infty}) \int_{D}(\a\circ\eta^{-1}_{ts})(x)|\nabla (g\circ\eta^{-1}_t)|\dx\\
	&\leq (1+|t|\,\|H\|_{\infty}) \int_{D}\a(\eta^{-1}_{ts}\circ\eta_{t})(x)|(D\eta_{t})^{-1}|\,|\nabla g|\,|\det D\eta_{t}|\dx\displaybreak[1]\\
	&\leq (1+|t|\,\|H\|_{\infty})^2 \Big(1+|t|\,\|\div w\|_{\infty} \\
	&\qquad\qquad\qquad\qquad\qquad +|t|\,\|H\|_{\infty}+O(t^2)\Big)\int_{D}\a(\eta^{-1}_{ts}\circ\eta_{t})(x)|\nabla g|\dx.
\end{align*}
Also, it can be checked that $ |\a(\eta^{-1}_{ts}\circ \eta_t)(x)-\a(x)|\lesssim [\a]_{1,K} |t| $ for $ x\in K $, since we have that $ |(\eta^{-1}_{ts}\circ \eta_t)(x)-x|\lesssim |t|(1+|s|)\|w\|_{L^{\infty}(K)} $. Hence,
\[
	\int_{D}\a(\eta^{-1}_{ts}\circ\eta_{t})(x)|\nabla g|\dx\leq \int_{D}\a(x)|\nabla g|\dx+C\,|t|\,[\a]_{1,K}\int_{D}|\nabla g|,
\]
where the constant $ C $ depends on $ D $ only. Recalling~\eqref{eq:lem1Rigot:2}-\eqref{eq:lem1Rigot:3}-\eqref{eq:lem1Rigot:4}, and the fact that $ |t_v|\leq C'|v| $, we can choose 
\[
	Q_{\a}=2(\sup_D\a)(\max\{1,C\}+[\a]_{1,D})\int_{D}|\nabla\chi_{L}| 
\] 
and this concludes the proof.
\end{proof}

\begin{remark}\label{rmrk:Rigot:2.1.1}
We note that in the proof of Lemma~\ref{lem:Rigot:2.1.1} above we only use the Lipschitzianity of $\a$.
\end{remark}

\medskip

Let us continue with an adaptation of a classical notion in geometric measure theory, to our setting with weight function $ \a $. Given $ x\in\R^d $ and $ r>0 $ let us define \emph{the weighted relative density function} of the set $ \Om_\gamma $ as
\[  
h_{\a}(x,r) \defeq \min\left\{\frac{\La{\Om_\gamma\cap B_r(x)}}{\La{B_r(x)}},\;\frac{\La{B_r(x)\setminus\Om_\gamma}}{\La{B_r(x)}}
\right\}.
\]
The rest of the proof is devoted to establishing a uniform lower bound of the form $ h_{\a}(x,r)\geq \e_0>0 $ for any point $ x\in\pt^*\Om_\gamma $, as long as $ r $ is taken sufficiently small. From here we will conclude the boundedness of the minimizer $ \Om_\gamma $ of $ \en $. 

Now, in view of the behavior of the density at infinity, $ \lim_{|x|\to\infty}\a(x)=+\infty $, the constant $ \Ra \ge 1 $ in condition~\eqref{eq:growth} can be chosen large enough so that $\R^d\setminus B_{\Ra}\subset \{x\in\Rd \colon \a(x)\geq 1\}$.
If $\Omega_\gamma\subset B_{2\Ra}$, then the minimizer is essentially bounded and we are done.  Therefore, in the following we may assume that the total variation measure $ |\nabla\chi_{\Om_\gamma}|=\mathscr{H}^{d-1}\mres \partial^*\Om_\gamma $ is nonvanishing in $\R^d\setminus B_{\Ra}$.  
We may also fix a constant $t_0\in (0,1)$ and balls $ B^{1} $ and $ B^{2} $, each of radius $ t_0 $, for which

\begin{equation}\label{eq:ball:cond}
	 3B^{1}\cap B^{2}=\emptyset,\quad 3B^{1}\cup B^{2}\subset\!\subset \R^d\setminus B_{\Ra} \quad\text{ and }\quad \int_{B^i}\a(x)|\nabla\chi_{\Om_\gamma}|>0\;\text{ for }\;i=1,2.	
\end{equation}

In what follows, $ B^1 $ and $ B^2 $ are to be used as reference sets, inside of which we will perform small deformations of our minimizer $ \Omega_\gamma $ in order to create a competitor set $ F $ with $ |F|=|\Om_\gamma| $ (see Lemma~\ref{lem:Rigot:2.1.1} above with $ D $ being $ B^1 $ or $ B^2 $, by analyzing two cases). This will allow us to exploit the volume constrained quasi-minimality of $ \Om_\gamma $ with respect to the $ \a $-perimeter, to derive a delicate growth estimate for the weighted relative density function of $ \Om_\gamma $ as a function of the radius $ r $, which will ultimately justify the uniform lower bound on $ h_{\a} $ that was claimed above. 

Before we continue, let us remark that for densities $ \a$ satisfying the assumptions \ref{itm:A1} and \ref{itm:A2} (hence, in particular, the condition \eqref{eq:growth}), the $ \a $-volume of any two sets $ F_1, F_2\subset\!\subset B_{2\bar{R}}\setminus B_{\bar{R}} $ are uniformly comparable:
\begin{equation}\label{eq:ball:comp}
	\Ca^{-1}\;\frac{|F_1|}{|F_2|}\leq \frac{\inf_{F_1}\a}{\sup_{F_2}\a}\,\frac{|F_1|}{|F_2|}\leq \frac{\La{F_1}}{\La{F_2}}\leq 	 \frac{\sup_{F_1}\a}{\inf_{F_2}\a}\,\frac{|F_1|}{|F_2|}
	\leq \Ca\;\frac{|F_1|}{|F_2|}
\end{equation}
for any $ \bar{R}>\Ra $. In particular, for any set $ F\subset\!\subset \R^d\setminus B_{\Ra} $,
\begin{equation}\label{eq:ball:comp:2}
	|F|=\int_{F}1\dx\leq\int_{F}\a(x)\dx= \La{F}
\end{equation}
These facts will be used in the following technical results. 

The first lemma establishes a bound on the growth rate of the weighted relative density function, for any minimizer $ \Om_\gamma $ of $ \en $, as a function of the radius $ r $ on balls having small $ \a $-volume, provided that the set $ \Om_\gamma $ or its complement $ \R^d\setminus\Om_\gamma $ have small density on that ball.

\blemma\label{lem:Rigot:2.1.2}
	For every $ \gamma>0 $ there exist $ 0<\e'<1 $ with $ \max\{1,\gamma^d\}\e' \ll 1 $, and $ 0<t_0<1 $, such that, for any minimizer $ \Om_\gamma $ of $ \en $ and  for any ball $ B_r(x)\subset\!\subset\R^d\setminus B_{R_{\a}}$ with $ 0<r<t_0 $ and $ \Ra $ as in condition~\eqref{eq:growth}, there holds: If $h_{\a}(x,r)<\e'\,\text{ and }\,0<\La{B_r(x)}\leq 1$, then 
	\[
	h_{\a}\left(x,\frac{r}{2}\right)\leq \frac{1}{2}h_{\a}(x,r).  
	\]
\elemma

\begin{proof}
Let $ x\in\R^d $ and $ 0<r<t_0 $. Loosely speaking, the main strategy is the following. If $ h_{\a}(x,r)=\La{\Om_\gamma\cap B_r(x)}/\La{B_r(x)} $ is small, we would like to delete the portion of $ \Om_\gamma $ that is inside of $ B_{t_*}(x) $ for some $ t_*\in (r/2,r) $ appropriately chosen. In a similar fashion, if $ h_{\a}(x,r)=\La{B_r(x)\setminus \Om_\gamma}/\La{B_r(x)} $ is small, we would like to append the ball $ B_{t_*}(x) $ to $ \Om_\gamma $. In these two cases, the resulting set has an additional portion of its boundary located inside $ \ol{B}_{t_*}(x) $, when compared to $ \Om_\gamma\cap B_r(x) $ or $ B_r(x)\setminus\Om_\gamma $. As $ h_{\a}(x,r) $ is assumed small, we lose in the volume term in~\eqref{eq:quasi-minimality} less than what we are adding on the boundary ($ F $ in~\eqref{eq:quasi-minimality} is our resulting set). We must then analyze the contribution of the boundary term in $ \ol{B}_{t_*}(x) $. We choose $ t_{*} $ in such a way that we can control $ \Ha{\Om_\gamma\cap\partial B_{t_*}(x)} $ in the first case, and $ \Ha{\partial B_{t_*}(x)\setminus\Om_\gamma} $ on the second case, in terms of $ \La{B_r(x)}^{(d-1)/d}\,h_{\a}(x,r) $. 
\smallskip

We distinguish four cases.
\smallskip

\noindent{\emph{Case 1:}} Assume that $ h_{\a}(x,r)=\La{\Om_\gamma\cap B_r(x)}/\La{B_r(x)} $ and $ B_r(x)\cap B^{1}=\emptyset $, where $ B^1 $ denotes the fixed ball in~\eqref{eq:ball:cond}. By Fubini and Chebyshev inequalities one can find $ t_*\in (r/2,r) $ and $ C=C(d)>0 $ such that
\begin{equation}\label{eq:lem2Rigot:1}  
	\Ha{\Om_\gamma\cap \partial B_{t_*}(x)}\leq C\,\frac{ (\sup_{B_r(x)}\a)^{1/d}}{\La{B_r(x)}^{1/d}}\;\La{\Om_\gamma\cap B_r(x)},
\end{equation}
Indeed, using Chebyshev's inequality with $ M_{\theta} \defeq \big(\theta\,\La{B_{r}(x)}^{1/d}\big)^{-1}\La{\Om_\gamma\cap B_r(x)} $, for $ \theta>0 $, yields 
\begin{align*}
	\Big|\{t\in (r/2,r) \colon \Ha{\Om_\gamma\cap \partial B_t(x)}>M_{\theta}\}\Big|
	&\leq \frac{1}{M_{\theta}}\int^r_{r/2}\Ha{\Om_\gamma\cap \partial B_t(x)} \d t \\
	&\leq \frac{1}{M_{\theta}}\La{\Om_\gamma\cap (B_{r}(x)\setminus B_{r/2}(x))} \\
	&\leq \theta\,\La{B_r(x)}^{1/d} \leq\theta\,\biggl(\sup_{B_r(x)}\a\biggr)^{1/d} \omega^{1/d}_d\,r< \frac{r}{2}
\end{align*} 
for $\theta \defeq (4^d\omega_d\sup_{B_r(x)}\a)^{-1/d}$.
On the other hand, note that
\[  
	\int_{\R^d}\a(y)|\nabla\chi_{\Om_\gamma\setminus B_{t_*}(x)}|=\int_{\R^d\setminus B_{t_*}(x)}\a(y)|\nabla\chi_{\Om_\gamma}|+\Ha{\Om_\gamma\cap\partial B_{t_*}(x)}.
\]
With $ D=B^{1} $, as defined in \eqref{eq:ball:cond}, and  $ L=\Om_\gamma $, we apply Lemma~\ref{lem:Rigot:2.1.1} to obtain values $ \e $ and $ Q_{\a} $ satisfying the conclusions of that lemma.
The constants $ \e$ and $ Q_{\a} $ depend on $ \Om_\gamma\cap B^{1} $ and $ B^{1} $ only, so, in particular, they are independent of $ x $ and $ r $ (and $ t_* $). 
Fix $ \e'\in(0,\e) $ for the time being, and choose $ r\in (0,t_0) $ sufficiently enough so that $ h_{\a}(x,r)<\e' $ and $ |B_r(x)|<\e $. Later on $ \e' $ will be reduced accordingly (independent of $ r $). 

Let $\Omega' \defeq \Om_\gamma\setminus B_{t_*}(x)$
and observe, by hypothesis, that $ \Omega'\cap B^{1}=\Om_\gamma\cap B^{1} $. This means that the same constants $ \e $ and $ Q_{\a} $ still work for Lemma~\ref{lem:Rigot:2.1.1} when applying it with $ L=\Omega' $ and $ D=B^{1} $. Applying Lemma~\ref{lem:Rigot:2.1.1} to $ L=\Omega' $, $ D=B^{1} $ and $ v=|\Om_\gamma\cap B_{t_*}(x)|$ (note that $ |v|<|B_r(x)|<\e $), we find a set $ F\subset\R^d $ such that $ F=\Omega' $ in a neighborhood of $ \R^d\setminus B^{1} $, $ |F|=|\Omega'|+|\Om_\gamma\cap B_{t_*}(x)|=1 $ and for which
\beqn \label{eq:aux1}
\begin{gathered}
\int_{B^{1}}\a(y)|\nabla\chi_{F}|
\leq\int_{B^{1}}\a(y)|\nabla\chi_{\Omega'}|+Q_{\a}|\Om_\gamma\cap B_{t_*}(x)| \\
 \text{ and } \quad
|F\triangle \Omega'| \leq Q_{\a}|\Om_\gamma\cap B_{t_*}(x) |,
\end{gathered}
\eeqn
where $ Q_{\a}\lesssim \big( \sup_{B^{1}}\a \big)\big(1+[\a]_{1,B^{1}}\big) $ is fixed and independent of $ x $ and $ r $.

In light of~\eqref{eq:ball:comp:2}, we have $|\Om_\gamma\cap B_{t_*}(x)|
\leq \La{\Om_\gamma\cap B_{t_*}(x)}$,
and in addition, using $ h_{\a}(x,r)<\e' $ and $ \La{B_r(x)}\leq 1 $, we deduce that
\begin{equation}\label{eq:lem2Rigot:2}
\La{\Om_\gamma\cap B_{t_*}(x)}
\leq \La{\Om_\gamma\cap B_{r}(x)}\\
= \La{B_r(x)}\;h_{\a}(x,r)\leq \La{B_r(x)}^{(d-1)/d}\;h_{\a}(x,r)< \e'.
\end{equation}  
Since $ F=\Omega' $ in a neighborhood of $ \R^d\setminus B^{1} $, we have
\[  
	\int_{\R^d\setminus B^{1}}\a(y)|\nabla\chi_F|=\int_{\R^d\setminus B^{1}}\a(y)|\nabla\chi_{\Omega'}|.
\]
Combining these facts, we obtain
\begin{align*} 
	\int_{\R^d}\a(y)|\nabla\chi_{F}|
	&\,\,\, \leq \int_{\R^d}\a(y)|\nabla\chi_{\Omega'}|+Q_{\a}\La{\Om_\gamma\cap B_{t_*}(x)}\displaybreak[1]\\
	&\xupref{eq:lem2Rigot:2}{\leq} \int_{\R^d}\a(y)|\nabla\chi_{\Omega'}|+Q_{\a}\La{B_r(x)}^{(d-1)/d}\,h_{\a}(x,r).
\end{align*}
In addition,
\begin{align*}
	\int_{\R^d}\a(y)|\nabla\chi_{\Omega'}|
	&\,\,\, =\int_{\R^d\setminus B_{t_*}(x)}\a(y)|\nabla\chi_{\Om_\gamma}|+\Ha{\Om_\gamma\cap\partial B_{t_*}(x)}\\
	&\xupref{eq:lem2Rigot:1}{\leq} \int_{\R^d\setminus B_{t_*}(x)}\a(y)|\nabla\chi_{\Om_\gamma}|+C\,\frac{(\sup_{B_r(x)}\a)^{1/d}}{\La{B_r(x)}^{1/d}}\;\La{\Om_\gamma\cap B_r(x)}\\
	&\xupref{eq:lem2Rigot:2}{\leq} \int_{\R^d\setminus B_{t_*}(x)}\a(y)|\nabla\chi_{\Om_\gamma}|+C\,(\sup_{B_r(x)}\a)^{1/d}\La{B_r(x)}^{(d-1)/d}\;h_{\a}(x,r).
\end{align*}
Assuming, without loss of generality, that $ \sup_{B_r(x)}\a \gg 1 $, the above estimate yields
\begin{equation}\label{eq:lem2Rigot:3}  
	\int_{\R^d}\a(y)|\nabla\chi_{F}|
	\leq \int_{\R^d\setminus B_{t_*}(x)}\a(y)|\nabla\chi_{\Om_\gamma}| +C \big(\sup_{B_r(x)}\a \big)^{1/d}\La{B_r(x)}^{(d-1)/d}\;h_{\a}(x,r)
\end{equation}
On the other hand, recalling~\eqref{eq:ball:comp:2}, \eqref{eq:aux1}, and~\eqref{eq:lem2Rigot:2}, we have
\begin{align*}
	|F\triangle \Om_\gamma| &\leq |F\triangle \Omega'|+|\Omega'\triangle\Om_\gamma| \leq (1+Q_{\a})|\Om_\gamma\cap B_{t_*}(x)| \\
	&\leq (1+Q_{\a})\La{\Om_\gamma\cap B_{t_*}(x)} \leq(1+Q_{\a})\La{B_r(x)}\,h_{\a}(x,r)\leq C\,\e'.
\end{align*}
Recall now that $ \Om_\gamma $ is a volume constrained quasi-minimal set for $ \a $-perimeter with $ g_{\gamma}(x)\simeq \gamma\,x=o(x^{(d-1)/d}) $; see~\eqref{eq:quasi-minimality}. By reducing the value of $ \e' $ in such a way that $ \max\{1,\gamma\}{\e'}^{\frac{1}{d}} \ll 1 $, we obtain
\[
	g_{\gamma}(|F\triangle \Om_\gamma|)\leq \eta \La{B_r(x)}^{\frac{d-1}{d}}h_{\a}(x,r)^{\frac{d-1}{d}},
\]
where $ \eta $ will be chosen below (independent of $ \e' $). Since $ \Om_\gamma $ is a volume constrained quasi-minimizer for the $ \a $-perimeter, the above inequality and~\eqref{eq:lem2Rigot:3} yield the following estimate:
\begin{equation}\label{eq:lem2Rigot:4}   
	\int_{B_{t_*}(x)}\a(y)|\nabla\chi_{\Om_\gamma}|\leq \biggl(C \big( \sup_{B_r(x)}\a \big)^{1/d}\,h_{\a}(x,r)+\eta\, h_{\a}(x,r)^{\frac{d-1}{d}}\biggr)\La{B_r(x)}^{\frac{d-1}{d}}.
\end{equation} 
In view of~\eqref{eq:lem2Rigot:4}, using the standard isoperimetric inequality on balls, we obtain that
\begin{align*}
	& \La{B_{r/2}(x)}\, h_{\a}\left(x,r/2\right) \\[0.5em]
	&\qquad\qquad =\min\left\{\La{\Om_\gamma\cap B_{r/2}(x)},\La{B_{r/2}(x)\setminus\Om_\gamma}\right\}\\[0.5em]
	&\qquad\qquad  \leq \sup_{B_{r/2}(x)}\a\cdot\min\left\{|\Om_\gamma\cap B_{r/2}(x)|,|B_{r/2}(x)\setminus\Om_\gamma|\right\}\\
	&\qquad\qquad  \xupref{eq:growth}{\leq}\big(\Ca\inf_{B_{r/2}(x)}\a\big)\cdot\; C\left(\int_{B_{r/2}(x)}|\nabla\chi_{\Om_\gamma}|\right)^{\frac{d}{d-1}}\\
	&\qquad\qquad  \leq C\,\Ca\big(\inf_{B_{r/2}(x)}\a\big)^{-\frac{1}{d-1}}\left(\int_{B_{t_*}(x)}\a(y)|\nabla\chi_{\Om_\gamma}|\right)^{\frac{d}{d-1}}\displaybreak[1]\\
	&\qquad\qquad  \leq C\,\Ca\;(\inf_{B_{r/2}(x)}\a)^{-\frac{1}{d-1}}\biggl(C\big(\sup_{B_r(x)}\a\big)^{1/d}\,h_{\a}(x,r)+\eta\, h_{\a}(x,r)^{\frac{d-1}{d}}\biggr)^{\frac{d}{d-1}}\La{B_r(x)}\\
	&\qquad\qquad  \leq C\,\Ca\;(\inf_{B_{r/2}(x)}\a)^{-\frac{1}{d-1}}\biggl(\big(\sup_{B_r(x)}\a\big)^{\frac{1}{d-1}}\,h_{\a}(x,r)^{\frac{1}{d-1}}+\eta^{\frac{d}{d-1}}\biggr)h_{\a}(x,r)\La{B_r(x)}\\
	&\qquad\qquad  \leq C\,\Ca\;\left(\frac{\sup_{B_{r}(x)}\a}{\inf_{B_{r/2}(x)}\a}\right)^{\frac{1}{d-1}}\biggl({\e'}^{\frac{1}{d-1}}+\eta^{\frac{d}{d-1}}\biggr)\;h_{\a}(x,r)\La{B_r(x)}\\[0.25em]
	&\qquad\qquad \xupref{eq:ball:comp}{\leq} C\,\Ca^2\;\biggl({\e'}^{\frac{1}{d-1}}+\eta^{\frac{d}{d-1}}\biggr)\;h_{\a}(x,r)\big(2^d\,\Ca\,\La{B_{r/2}(x)}\big).
\end{align*}
In other words, we conclude that
\[  
	h_{\a}\left(x,\frac{r}{2}\right) \leq C\;\Ca^3\;\Big({\e'}^{\frac{1}{d-1}}+\eta^{\frac{d}{d-1}}\Big)\;h_{\a}(x,r) \leq \frac{1}{2}h_{\a}(x,r),
\]
provided $ 0<\eta<(4C\Ca^3)^{-\frac{d-1}{d}} $ is chosen small enough; and possibly reducing the value of $ \e'\lesssim \min\{\e,(4C\Ca^3)^{-(d-1)}\} $, independently of $ x $ and $ r $.

\bigskip

\noindent{\emph{Case 2:}} Suppose $ h_{\a}(x,r)=\La{\Om_\gamma\cap B_r(x)}/\La{B_r(x)} $ and $ B_r(x)\cap B^{1}\neq\emptyset $. Under the hypotheses, $ B_{r}(x)\subset 3B^{1} $, and~\eqref{eq:ball:cond} ensures that $ B_r(x)\cap B^{2}=\emptyset $. We proceed exactly as in the argument in the first case, with $ B^{1} $ replaced by $ B^{2} $.

\bigskip

\noindent{\emph{Case 3:}} Assume that $ h_{\a}(x,r)=\La{B_r(x)\setminus \Om_\gamma}/\La{B_r(x)} $ and $ B_r(x)\cap B^{1}=\emptyset $. We will proceed analogously as the first case, but with 
\[ 
	\Omega' \defeq \Om_\gamma\cup B_{t_*}(x) 
\]
for a certain $ t_{*} $ to be determined below. Again, by Fubini and Chebyshev inequalities, we can find $ t_{*}\in (r/2,r) $ in such a way that
\begin{equation}\label{eq:lem2Rigot:5}  
	\Ha{\partial B_{t_*}(x)\setminus\Om_\gamma}
	\leq C\, \frac{(\sup_{B_r(x)}\a)^{1/d}}{\La{B_r(x)}^{1/d}}\La{B_r(x)\setminus\Omega_\gamma}.
\end{equation} 
Also, we have that
\begin{equation}\label{eq:lem2Rigot:6}  
	\int_{\R^d}\a(y)|\nabla\chi_{\Om_\gamma\cup B_{t_*}(x)}|=\int_{\R^d\setminus B_{t_*}(x)}\a(y)|\nabla\chi_{\Om_\gamma}|+\Ha{\partial B_{t_*}(x)\setminus\Om_\gamma}.
\end{equation}
Fix now $ \e' $ as in the previous cases so that $ h_{\a}(x,r)<\e'<\e $ and $ |B_r(x)|<\e $. Note that~\eqref{eq:ball:comp:2} yields $|\Omega'\setminus\Om_\gamma|\leq\La{\Omega'\setminus \Om_\gamma}$. Also, recalling that $ \La{B_r(x)}\leq 1 $, we have 
\begin{align*}
	\La{\Omega'\setminus\Om_\gamma}
	&\leq\La{B_{t_*}(x)\setminus\Om_\gamma}\leq \La{B_{r}(x)\setminus\Om_\gamma}\\
	&\leq\La{B_r(x)}\;h_{\a}(x,r)\leq \La{B_r(x)}^{(d-1)/d}\;h_{\a}(x,r)< \e'.
\end{align*}
Because of the choice of $ \e' $, we can apply Lemma~\ref{lem:Rigot:2.1.1} to $ L=\Omega' $, $ D=B^{1} $ and $ v=-|\Omega'\setminus\Om_\gamma| $ (with $ |v|\leq |B_r(x)|<\e $) to obtain the existence of a set $ F\subset\R^d $ such that $ F=\Omega' $ in a neighborhood of $ \R^d\setminus B^{1} $ and $ |F|=|\Omega'|-|\Omega'\setminus\Om_\gamma|=1 $. Moreover, $F$ satisfies the estimates
\[  
\begin{split}
	\int_{B^{1}}\a(y)|\nabla\chi_F|
	&\leq\int_{B^{1}}\a(y)|\nabla\chi_{\Omega'}|+Q_{\a}|\Omega'\setminus\Om_\gamma|\\
	&\leq\int_{B^{1}}\a(y)|\nabla\chi_{\Omega'}|+Q_{\a}\,\La{B_r(x)}^{(d-1)/d}\;h_{\a}(x,r),
\end{split}
\]
and
\[
	|F\triangle\Omega'|
	\leq Q_{\a}\,|\Omega'\setminus\Om_\gamma|,
\]
where $ Q_{\a}\lesssim \big(\sup_{B^{1}}\a \big) \big(1+[\a]_{1,B^{1}}\big) $ is fixed, and independent of $ x $ and $ r $.

Then, just like in the first case, we deduce
\begin{equation}\label{eq:lem2Rigot:7}   
	\int_{\R^d}\a(y)|\nabla\chi_F|\leq \int_{\R^d}\a(y)|\nabla\chi_{\Omega'}|
	+C\La{B_r(x)}^{(d-1)/d}\;h_{\a}(x,r).
\end{equation}
Thus, utilizing~\eqref{eq:lem2Rigot:6} and \eqref{eq:lem2Rigot:5}, we can estimate
\[  
\begin{split}
	\int_{\R^d}\a(y)|\nabla\chi_{\Omega'}|
	&\leq\int_{\R^d\setminus B_{t_*}(x)}\a(y)|\nabla\chi_{\Om_\gamma}|+
	C\,\frac{(\sup_{B_r(x)}\a)^{1/d}}{\La{B_r(x)}^{1/d}}\La{B_r(x)\setminus\Omega_\gamma}\\
	&\leq\int_{\R^d\setminus B_{t_*}(x)}\a(y)|\nabla\chi_{\Om_\gamma}|+
	C\,(\sup_{B_r(x)}\a)^{1/d}\La{B_r(x)}^{(d-1)/d}\,h_{\a}(x,r),
\end{split}	
\]
which can be combined with~\eqref{eq:lem2Rigot:7}, in turn, to obtain
\[  
	\int_{\R^d}\a(y)|\nabla\chi_F|\leq \int_{\R^d\setminus B_{t_*}(x)}\a(y)|\nabla\chi_{\Om_\gamma}|
	+C\,(\sup_{B_r(x)}\a)^{1/d}\,\La{B_r(x)}^{(d-1)/d}\;h_{\a}(x,r).
\]
On the other hand, 
\[  
\begin{split}
	|F\triangle \Om_\gamma|
	&\leq |F\triangle \Omega'|+|\Omega' \triangle \Om_\gamma|\\
	&\leq (Q_{\a}+1)|\Omega'\setminus\Om_\gamma|\leq(Q_{\a}+1)\La{\Omega'\setminus\Om_\gamma}\\
	&\leq (Q_{\a}+1)\,\La{B_r(x)}\;h_{\a}(x,r)\leq C\,\e'.
\end{split}	
\] 
We reach the conclusion like in the first case, utilizing the volume constrained quasi-minimality of $ \Om_\gamma $, the relative isoperimetric inequality on balls, and the fact that the weight has uniformly bounded oscillation (condition~\eqref{eq:growth}).

\bigskip

\noindent{\emph{Case 4:}} Suppose $ h_{\a}(x,r)=\La{B_r(x)\setminus \Om_\gamma}/\La{B_r(x)} $ and $ B_r(x)\cap B^{1}\neq\emptyset $. We employ the same construction as in the third case, with $ B^{1} $ replaced by $ B^{2} $, to conclude the proof of this lemma.
\end{proof}

\bigskip

\blemma\label{lem:Rigot:2.1.3}
	Given $ \gamma>0 $, there exist $ 0<\e_0<1 $ with $ \max\{1,\gamma^d\}\e_0 \ll 1 $, and $ 0<t_0<1 $, such that for any minimizer $ \Om_\gamma $ of $ \en $ and any ball $
	 B_{r}(x)\subset\!\subset\R^d\setminus B_{\Ra} $, for which $ 0<r\leq t_0 $, and for $ \Ra $ as in~\eqref{eq:growth}, the following holds: If $ h_{\a}(x,r)<\e_0\quad\text{ and }\quad 0<\La{B_{r}(x)}\leq 1$, then 
	\[
		\text{ either }\;\; \left|\Om_\gamma\cap B_{r/2}(x)\right|=0\;\;\text{ or }\;\; \left|B_{r/2}(x)\setminus\Om_\gamma\right|=0.
	\]
\elemma

\medskip

\begin{remark}\label{rmrk:Rigot:2.1.3}
We note that, in particular, for any point $ x\in\pt^*\Om_\gamma $ and any ball $ B_{r}(x) $ as above, one has
	\[  
		h_{\a}(x,r)=\frac{\La{\Om_\gamma\cap B_r(x)}}{\La{B_r(x)}}\geq \e_0>0.
	\]
\end{remark}

\medskip

\begin{proof}
Under the hypotheses of this lemma, for any $ y\in B_{r/2}(x) $ we write $ d_k \defeq h_{\a}(y,2^{-k}r) $. Choose 
\[
	\e_0 \defeq (2^d\Ca)^{-1}\min\{\e',1/2\},
\]
with $ \e'>0 $ given as in  Lemma~\ref{lem:Rigot:2.1.2}. If $ h_{\a}(x,r)=\La{B_r(x)}^{-1}\;\La{\Om_\gamma\cap B_{r}(x)} $, then
\[  
	d_1 \defeq {\La{B_{r/2}(x)}}^{-1}\La{\Om_\gamma\cap B_{r/2}(y)}\leq 2^d\Ca\;h_{\a}(x,r)<\min\{\e',1/2\}
\]
where we used the estimate~\eqref{eq:ball:comp}. Since $ d_1<\e' $, by Lemma~\ref{lem:Rigot:2.1.2} we deduce that
\[  
	d_2=h_{\a}\left(y,\frac{r}{4}\right) \leq \frac{1}{2}h_{\a}\left(y,\frac{r}{2}\right)=\frac{d_1}{2}.
\]
In particular, $ d_2<\min\{\e',1/2\}$. Hence, by induction it follows that, for all $k\geq 1$, we have $d_{k+1}\leq d_1/2^k$.

Hence, for any $ y\in B_{r/2}(x) $, and by virtue of condition~\eqref{eq:growth}, for all $k \geq 1$, we have
\[ 
\frac{\e'}{2^{k-1}}
	\geq \frac{\La{\Om_\gamma\cap B_{2^{-k}r}(y)}}{\La{B_{2^{-k}r}(y)}}
	\geq\frac{\inf_{B_r(y)}\a}{\sup_{B_r(y)}\a} \, \frac{|\Om_\gamma\cap B_{2^{-k}r}(y)|}{|B_{2^{-k}r}(y)|}\geq\Ca^{-1}\,\frac{|\Om_\gamma\cap B_{2^{-k}r}(y)|}{|B_{2^{-k}r}(y)|}
\]
Thus, $ y $ is not a point of density of $ \Om_\gamma $, which yields $ |\Om_\gamma\cap B_{r/2}(x)|=0 $.

Otherwise, if $ h_{\a}(x,r)=r^{-d}\La{B_{r}(x)\setminus \Om_\gamma} $, it can be shown in the same fashion that, for any $ k\geq 1 $ and any $ y\in B_{r/2}(x) $,
\[  
	d_{k}=\La{B_{2^{-k}r}(y)}^{-1}\La{B_{r}(y)\setminus\Om_\gamma}<2^d\,\Ca\;\e<\min\{\e', 1/2\}\quad\text{ and }\quad d_{k+1}\leq \frac{d_k}{2}.
\]
Therefore,
\[ 
\Ca^{-1}\cdot\frac{|B_{2^{-k}r}(y)\setminus\Om_\gamma |}{|B_{2^{-k}r}(y)|}\leq  \frac{\La{B_{2^{-k}r}(y)\setminus\Om_\gamma}}{\La{B_{2^{-k}r}(y)}}
\leq\frac{\e'}{2^{k-1}}\underset{k\to\infty}{\longrightarrow} 0.
\]
and so $ y $ is not a density point of $ \R^d\setminus\Om_\gamma $. This shows that $ |B_{r/2}(x)\setminus\Om_\gamma|=0 $.
\end{proof}

\bigskip

\begin{proof}[Proof of Theorem~\ref{thm:bd} (Part 1)]

We now finish the argument to establish the essential boundedness of the minimizer $ \Om_\gamma $ in Theorem~\ref{thm:bd} under the hypotheses \ref{itm:A1} and \ref{itm:A2}. Let $ \e_0 $ and $ t_0>0 $ be given as in Lemma~\ref{lem:Rigot:2.1.3}, and let $j_{\a}\in\mathbb{N}$ be such that $ \Ra<2^{j_{\a}}$.

Let $ \{R_1,\ldots, R_N\}\subset (2^{j_{\a}},+\infty) $ be any finite collection for which 
\[
	\pt^*\Om_\gamma\cap (B_{R_{i+1}}\setminus\overline{B}_{R_i})\neq \emptyset\;\;\text{ and }\;\; R_{i+1}-R_i\geq 1,\quad \text{ for all }i.
\] 
Then we can select balls $ B_{r_i}(x^i)\subset\!\subset B_{R_{i+1}}\setminus B_{R_i} $ with $ x^i\in\pt^*\Om_\gamma $ and $ r_i\leq t_0< 1 $ taken suffieciently small so that 
\[  
	\frac{1}{2}\leq\La{B_{r_i}(x^i)}< 1. 
\]
In view of~\eqref{eq:pf:thm1:6}, and employing Remark~\ref{rmrk:Rigot:2.1.3} over each ball, we derive the following bound 
\[  
\begin{split}
	c_d\,\Ca\int_{\R^d\setminus B_{2^{j_{\a}}}}\a(x)|\nabla\chi_{\Om_\gamma}|
	&\geq \La{\Om_\gamma\setminus B_{2^{j_{\a}}}}\geq \sum^N_{i=1}\La{\Om_\gamma\cap (B_{R_{i+1}}\setminus B_{R_i})}\\
	&\geq \sum^N_{i=1}\La{\Om_\gamma\cap B_{r_i}(x^i)}\geq \sum^N_{i=1}\e_0\,\La{B_{r_i}(x^i)}\\
	&\geq \tfrac{1}{2}\e_0\, N
\end{split}
\]
thus proving that the maximal number of such radii $ \{R_i\} $ must be finite. This shows that there exists $ R_{\rm max}\in (\Ra,+\infty) $, depending on $ \a $, for which $ \pt^*\Om_\gamma\cap (\R^d\setminus B_{R_{\rm max}})=\emptyset $. Since $\Om_\gamma$ is assumed to be essentially closed, then
\[ 
	\Om_\gamma\subset B_{R_{\max}}.
\]
This establishes Theorem~\ref{thm:bd} under the hypotheses \ref{itm:A1} and \ref{itm:A2}.
\end{proof}

\medskip

\subsection{Radial and monotone densities}

Now we turn to densities $\a$ satisfying the assumptions \ref{itm:A1} and \ref{itm:A3}. A key ingredient in the second part of the proof of Theorem~\ref{thm:bd} is the $\e-\e^{(d-1)/d}$ property recently proved by Pratelli and Saracco \cite{PrSa}. We state this result as a lemma here for the convenience of the reader.

\blemma[cf. Theorem A in \cite{PrSa}] \label{lem:eps_epsbeta} For any set of finite perimeter $E\subset \Rd$ and for any ball $B$ such that $\Hn(B \cap \pt^*E)>0$ there exists $\ol{\e}>0$ and $C_E>0$ such that for any $|\e|<\ol{\e}$ there is a set $F\subset\Rd$ satisfying
	\beqn \label{eq:eps_epsbeta}
		F \triangle E \subset\!\subset B, \qquad |F|-|E| = \e, \qquad \aPer{F}-\aPer{E} \leq C_E|\e|^{(d-1)/d}.
	\eeqn
Moreover the constant $C_E$ can be chosen arbitrarily small up to possibly choosing $\ol{\e}$ smaller.
\elemma

Now we will finish the proof of Theorem~\ref{thm:bd} following an argument similar to \cite[Theorem 5.9]{MoPr} and \cite[Lemma 5.1]{KnMu2014}.

\bigskip

\begin{proof}[Proof of Theorem~\ref{thm:bd} (Part 2)]
Assume, for a contradiction, that $\Om_\gamma$ is not bounded. Then $|\Om_\gamma \setminus B_r|>0$ for all $r>0$. In particular, for $r \geq R_\a$,
	\[
		|\Om_\gamma \setminus B_r|_{\a} \geq \a(R_{\a})|\Om_\gamma \setminus B_r| >0.
	\]
Now let
	\[
		\Om(r) \defeq \Om_\gamma \cap B_r, \quad \Om_r \defeq \Om_\gamma \cap \pt B_r, \quad
		\aPer{r} \defeq \Ha{\pt^* \Om_\gamma \setminus B_r}, \text{ and } \ V_\a(r) \defeq |\Om_\gamma \setminus B_r|_\a.
	\]
By the proof of Theorem~4.3 of \cite{MoPr}, 
$$\aPer{\Om(r)} < \aPer{\Om_\gamma} \quad\text{ for $r \geq R_\a$.}
$$
This follows from a nifty argument involving projection onto the sphere $B_r$ and for the reader’s convenience we include the details here. Consider the projection map $\Pi \colon \pt^*\Om_\gamma \setminus B_r \to \pt B_r$. Clearly, $\Pi$ is strictly 1-Lipschitz and its image $\img(\Pi)$ satisfies
	\[
		\pt^* \Om(r) \setminus \pt^*\Om_\gamma \subset \img(\Pi).
	\]
This inclusion would be trivially true for bounded $\Om_\gamma$; however, since we assumed that $\Om_\gamma$ is unbounded it contains the whole cone $C=\{\lambda x \colon \lambda \geq 1,\ x\in H\}$ where $H=\left(\pt^*\Om(r)\setminus\pt^*\Om_\gamma\right)\setminus \img(\Pi)$. But since $\a$ is increasing for large values of $|x|$, the cone $C$, and hence the set $\Om_\gamma$, have infinite volume, unless $\Hn(H)=0$. Thus $H=\emptyset$ up to $\Hn$-measure zero.

Since $\a$ is eventually increasing, for $r \geq R_\a$, by the co-area formula we get that
	\begin{align*}
		\aPer{\Om(r)} &= \int_{\pt^*\Om_\gamma \cap B_r} \a(|x|) \d\Hn(x) + \int_{\pt^*\Om(r) \setminus \pt^*\Om_\gamma} \a(\theta)\,d\Hn(\theta) \\
							  &< \int_{\pt^*\Om_\gamma \cap B_r} \a(|x|) \d\Hn(x) + \int_{\pt^*\Om_\gamma \setminus \pt B_r} \a(\Pi(x)) \d\Hn(x) \\
							  &\leq \int_{\pt^*\Om_\gamma} \a(|x|)\d\Hn(x) = \aPer{\Om_\gamma}.
	\end{align*}

Also, note that
	\[
		\aPer{\Om(r)} = \aPer{\Om_\gamma} - \aPer{r} + \Ha{\Om_r}.
	\]
This, combined with $\aPer{\Om(r)} < \aPer{\Om_\gamma}$, implies that
	\beqn \label{eq:bd-1}
		\aPer{r} > \Ha{\Om_r}.
	\eeqn

For $r \geq 1$, the standard isoperimetric inequality in the sphere yields that for any subset $E_r$ of the sphere $\pt B_r$ having area at most half of the sphere, we have $\Hnm(\pt^* E_r) \geq c \big( \Hn(E_r) \big)^{(d-2)/(d-1)}$ for some constant $c>0$ where $\pt^* E_r$ denotes the reduced boundary of $E_r$ inside $\pt B_r$. In fact, as $\Om_\gamma$ has finite volume, for $r$ sufficiently large, $\Hn(\Om_r) \leq \frac{1}{2} \Hn(\pt B_r)$. Since the density $\a$ is constant on $\pt B_r$, we get
	\beqn \label{eq:bd-2}
		\Hnm_\a(\pt^* \Om_r) \geq c \Big(\Ha{\Om_r}\Big)^{\frac{d-2}{d-1}} \big(\a(r)\big)^{\frac{1}{d-1}}.
	\eeqn
Combining \eqref{eq:bd-1} and \eqref{eq:bd-2}, we obtain
	\[
		\big(\aPer{r}\big)^{\frac{1}{d-1}}\, \Hnm_\a(\pt^* \Om_r) > \big(\Ha{\Om_r}\big)^{\frac{1}{d-1}}\, \Hnm_\a(\pt^*\Om_r) \geq c\, \Ha{\Om_r} \big(\a(r)\big)^{\frac{1}{d-1}}.
	\]
Since $\a(r) \geq \a(R_\a)$ for $r \geq R_\a$, the estimate above becomes
	\beqn \label{eq:bd-3}
		\Hnm_\a(\pt^* \Om_r) \geq c\, \big( \aPer{r} \big)^{-\frac{1}{d-1}} \, \Ha{\Om_r}.
	\eeqn
	
Also, observe that
	\beqn \label{eq:bd-4}
		-\frac{d}{dr}\aPer{r} = \left|\frac{d}{dr}\aPer{r}\right| \geq \Hnm_\a(\pt^*\Om_r) \quad \text{ and } \quad -\frac{d}{dr}V_\a(r) = \Ha{\Om_r}.
	\eeqn
Hence, by \eqref{eq:bd-3}, we get
	\[
		-\frac{d}{dr}\Big( \big(\aPer{r}\big)^{\frac{d}{d-1}} \Big)  \geq - c\, \frac{d}{dr}V_\a(r).
	\]
Since both $\aPer{r}$, $V_\a(r) \to 0$ as $r\to \infty$, integrating both sides from $r$ to $\infty$ yields
	\beqn \label{eq:bd-5}
		\big( \aPer{r} \big)^{\frac{d}{d-1}} \geq c\, V_\a(r).
	\eeqn
	
Now, let $R\in\R$ be such that $\Om_\gamma \cap B_R \neq \emptyset$. Then by Lemma~\ref{lem:eps_epsbeta}, there exists $\ol{\e}$ and $C_{\Om_\gamma}>0$ such that for all $|\e|<\ol{\e}$ there exists $\Om_\e$ satisfying
	\beqn \label{eq:bd-eps_epsbeta}
		\Om_\e \triangle \Om_\gamma \subset\!\subset B_R, \qquad |\Om_\e|-|\Om_\gamma| = \e, \qquad \aPer{\Om_\e}-\aPer{\Om_\gamma} \leq C_{\Om_\gamma}|\e|^{(d-1)/d}.
	\eeqn
By choosing $\e$ smaller, if necessary, by Lemma~\ref{lem:Rigot:2.1.1}  and Remark~\ref{rmrk:Rigot:2.1.1}, we can take $\Om_\e$ such that
	\beqn \label{eq:bd-giusti}
		|\Om_\e \triangle \Om_\gamma| < C\e.
	\eeqn
	
Let $r>R$ be such that $\e \defeq V_\a(r) < \ol{\e}$ and define $\Om_\e(r) = \Om_\e \cap B_r$. Then $|\Om_\e(r)|=|\Om_\gamma|=1$ and $\Om_\e(r)$ is competitor in \eqref{eq:nlip}. So,
	\begin{align*}
		\aPer{\Om_\e(r)} &= \aPer{\Om_\e} - \aPer{r} + \Ha{\Om_r} \\
								  &\upupref{eq:bd-5}{eq:bd-eps_epsbeta}{\leq} \aPer{\Om_\gamma} + C_{\Om_\gamma} \e^{\frac{d-1}{d}} - c\,\e^{\frac{d-1}{d}} + \Ha{\Om_r} \\
								  &\leq \aPer{\Om_\e(r)} + \gamma\, \bigg(\V(\Om_\e(r)) - \V(\Om_\gamma) \bigg) + C_{\Om_\gamma} \e^{\frac{d-1}{d}} - c\,\e^{\frac{d-1}{d}} + \Ha{\Om_r} \\
								  &\xupref{eq:V_Lip}{\leq} \aPer{\Om_\e(r)} + c \gamma \, |\Om_\e(r) \triangle \Om_\gamma| + C_{\Om_\gamma} \e^{\frac{d-1}{d}} - c\,\e^{\frac{d-1}{d}} + \Ha{\Om_r}.
	\end{align*}
Since, as noted in Lemma~\ref{lem:eps_epsbeta}, we can choose $C_{\Om_\gamma}$ arbitrarily small, and since \eqref{eq:bd-eps_epsbeta} and \eqref{eq:bd-giusti} hold, for $\e$ sufficiently small, we have
	\[
		\Ha{\Om_r} \geq C \e^{\frac{d-1}{d}}
	\]
for some $C>0$. Then, by \eqref{eq:bd-4}, $-\frac{d}{dr}V_\a(r) \geq C\e^{(d-1)/d} = C\big(V_\a(r)\big)^{(d-1)/d}$, i.e.,
	\[
		\frac{d}{dr} \left( \big(V_\a(r)\big)^{1/d} \right) \leq -C.
	\]
This contradicts the fact that $V_\a(r)>0$ for all $r>R_\a$ and we establish Theorem~\ref{thm:bd} under the hypotheses \ref{itm:A1} and \ref{itm:A3}.
\end{proof}

\addtocontents{toc}{\protect\setcounter{tocdepth}{2}}

%%%%%%%%%%%%%%%%%%%%%%%%%%%%%%%%%%%%%%%%%%%%%%%%%%%%%%%%%%%%%%%%%
%%%%%%%%%%%%%%%%%%%%%%%%%%%%%%%%%%%%%%%%%%%%%%%%%%%%%%%%%%%%%%%%%

\section{Global minimizers in the small $\gamma$ regime}\label{sec:globalmin}

In this section we will prove Theorem~\ref{thm:global}. The proof of this theorem relies on the regularity of almost-minimizers of the perimeter functional. As in \cite[Definition 2.5]{AcFuMo13}, we will call a set in $E\subset \Rd$  an \emph{$\omega$-minimizer} of the Euclidean perimeter functional $\per$ in an open set $O\subset\Rd$ with $\omega>0$ if for every ball $B_r(x)\subset O$ and for every set of finite perimeter $F \subset \Rd$ such that $E \triangle F \subset\!\subset B_r(x)$ we have
	\[
		\per(E) \leq \per(F) + \omega r^d.
	\]
This definition also coincides with the definition of \emph{almost-minimizers} according to \cite[Section 1.5]{tamanini} with $\alpha(r)=\omega r$, although it is weaker than the $(\omega,r)$-minimizers as defined in \cite{maggiBOOK}. Nevertheless, by the regularity theory developed in \cite{tamanini}, it will be sufficient to justify passage to a limit as $\gamma\to 0$.

First we show that any minimizer of the energy $\varE_\gamma$ is an $\omega$-minimizer of the Euclidean perimeter functional $\per$ in any open set that does not contain the origin. 

\blemma \label{lem:quasimin}
Let $\Om_\gamma\subset\Rd$ be a solution of the problem \eqref{eq:nlip}. Then $\Om_\gamma$ is $\omega$-minimal for the Euclidean perimeter functional $\per$ in $O=\Rd\setminus \ol{B_\delta}$ for any $\delta>0$.
\elemma

\begin{proof}
We first note that the constraint $|\Om_\gamma|=1$ may be replaced by a penalization term in the energy. Namely, arguing as in Step 1 of the proof of \cite[Theorem 2.7]{BoCr14} we obtain that there exists a constant $\lambda>0$ such that
	\[
		\min \calF_\gamma^\lambda = \calF_\gamma^\lambda(\Omega_\gamma) = \varE_\gamma(\Om_\gamma)
	\]
where
	\[
		\calF_\gamma^\lambda(E) = \varE_\gamma(E) + \lambda \big| |E|-1 \big|
	\]
for any $E\subset \Rd$.

Fix $r>0$ such that $B_r(x) \cap B_\delta = \emptyset$ and let $F \subset \Rd$ be any set such that $\Om_\gamma \triangle F \subset\!\subset B_r(x)$. Since $\varE_\gamma(\Om_\gamma) = \calF_\gamma^\lambda(\Om_\gamma) \leq \calF_\gamma^\lambda(F)$, using \eqref{eq:V_Lip} we obtain that
	\[
		\aPer{\Om_\gamma} \leq \aPer{F} + \omega |\Om_\gamma \triangle F|
	\]
for some constant $\omega$ depending only on $\gamma$ and $\lambda$.

Let $m \defeq \max \{ \a(y) \colon y\in B_r(x) \}$. Since $\a$ is Lipschitz there exists a constant $C>0$ such that $\a(y) \geq m-Cr$ for all $y\in B_r(x)$. Then
	\begin{align*}
		(m-Cr) \per(\Om_\gamma,B_r(x)) &\leq \aPer{\Om_\gamma,B_r(x)} \leq \aPer{F,B_r(x)} + \omega |\Om_\gamma \triangle F| \\
														   &\leq m \per(F,B_r(x)) + \omega |\Om_\gamma \triangle F|.
	\end{align*}
Hence,
	\[
		\per(\Om_\gamma,B_r(x)) \leq \per(F,B_r(x)) + Cr^d
	\]
provided $\per(\Om_\gamma,B_r(x)) \leq C r^{d-1}$. This establishes that $\Om_\gamma$ is an $\omega$-minimizer of the Euclidean perimeter in $\Rd\setminus \ol{B_\delta}$.

\smallskip

It now remains to prove that $\per(\Om_\gamma,B_r(x)) \leq C r^{d-1}$. Note that $\Om_\gamma$ satisfies the $\e-\e^{(d-1)/d}$ property given by \eqref{eq:eps_epsbeta} in Lemma~\ref{lem:eps_epsbeta}. Let $B^1$ and $B^2$ be two disjoint balls, intersecting $\pt^* \Om_\gamma$ in a set of positive $\Hn$-measure. Let  $C_1,\,C_2>0$ and $0<\e_1,\e_2<1$ be the corresponding constants satisfying the conditions in \eqref{eq:eps_epsbeta}. Take $C=\max\{C_1,C_2\}$, $\ol{\e} = \min\{\e_1,\e_2\}$, and $\ol{r} = \min \big\{ (\ol{\e}/\omega_d)^{1/d}, \dist(B^1,B^2) \big\}$.

If $r < \ol{r}$, then $\e \defeq |B_r(x) \cap \Om_\gamma| < \omega_d r^d \leq \ol{\e}$, and $B_r(x)$ cannot intersect both $B^1$ and $B^2$. Without loss of generality, assume that $B_r(x) \cap B^1 = \emptyset$. For this $\e$, let $F\subset \Rd$ be the set of finite perimeter satisfying \eqref{eq:eps_epsbeta}. Clearly, we also have
	\[
		\aPer{F} \leq \aPer{\Om_\gamma} + C\e \leq \aPer{\Om_\gamma} + C\omega_d r^d.
	\]
	
Now, let $G=F \setminus B_r(x)$. Then $|G|=|\Om_\gamma|=1$ and $G$ is admissible for \eqref{eq:nlip}. By minimality of $\Om_\gamma$ and Lipschitzianity of $\calV$ (see \eqref{eq:V_Lip}) we obtain
\begin{align*}
	\aPer{\Om_\gamma} &\leq \aPer{G} + \gamma \big( \N{G} - \N{\Om_\gamma} \big) \\
				 &\leq \aPer{F} - \aPer{\Om_\gamma,B_r(x)} + d\omega_d R_{\max}^p r^{d-1} + C\gamma \, |G \triangle \Om_\gamma| \\
				 &\leq \aPer{\Om_\gamma} + C\omega_d r^d  - \aPer{\Om_\gamma,B_r(x)} + d\omega_d R_{\max}^p r^{d-1} + C\gamma \, r^d \\
				 &\leq \aPer{\Om_\gamma} - \aPer{\Om_\gamma,B_r(x)} + C\,r^{d-1}
\end{align*}	 
with a constant depending only on $d$, $R_{\max}$, and $\gamma$, where $R_{\max}$ is the constant obtained in the proof of Theorem \ref{thm:bd}. Recalling that $\a(x)\geq \delta^p$ in $\Rd\setminus \ol{B_\delta}$ we have that
	\[
		\per(\Om_\gamma,B_r(x)) \leq  \delta^{-p} \aPer{\Om_\gamma,B_r(x)} \leq C r^{d-1}.
	\]
	
If $r \geq \ol{r}$, on the other hand, then
	\[
		\per(\Om_\gamma,B_r(x)) \leq \delta^{-p} \aPer{\Om_\gamma,B_r(x)} \leq \delta^{-p} \aPer{\Om_\gamma} \leq \frac{\delta^{-p}\aPer{\Om_\gamma}}{\ol{r}^{d-1}} \, r^{d-1}.
	\]
	
Hence, in either case we obtain that $\per(\Om_\gamma,B_r(x)) \leq C r^{d-1}$, completing the proof of the lemma.
\end{proof}

\medskip

Since the density $\a(x)=|x|^p$ is vanishing at zero we cannot conclude that $\Om_\gamma$ is an $\omega$-minimizer of the perimeter in $\Rd$. Nevertheless, combining \cite[Sections 1.9 and 1.10]{tamanini} (see also \cite{Almgren76,whiteminimax}), we obtain the following regularity result for any $\omega$-minimizer of perimeter in open sets.

\blemma \label{lem:reg} Let $O\subset\R^d$ be an open set.
\begin{enumerate}
\item  If $E\subset\R^d$ is an $\omega$-minimizer of Euclidean perimeter in $O$ then $\partial^*E\cap O$ is a $C^{1,\alpha}$ hypersurface for any $\alpha\in(0,1/2)$.
\item  If $E_n\subset\Rd$ is a sequence of uniformly $\omega$-minimizers of Euclidean perimeter in $O$ and if $E_n \to E$ in $L^1(O)$, then $\pt E_n \to \pt E$ in $C^{1,\alpha}$ and $\pt E_n$ is a $C^{1,\alpha}$ graph over $\pt E$.
\end{enumerate}
\elemma

\bigskip

We are now ready to present the proof of global minimality of balls in the small $\gamma$ regime.

\medskip

\begin{proof}[Proof of Theorem~\ref{thm:global}]
We fix $r_0=(1/\omega_d)^{1/d}$ with $|B_{r_0}|=1$, and recall (see e.g. \cite{AlBrChMePo,BrDePhilRu}) that $B_{r_0}$ is the unique minimizer of the local weighted isoperimetric problem $e(0)$.  
By Theorem~\ref{theorem:1}, for any $\gamma>0$, there exists $\Om_\gamma\subset \Rd$, $|\Om_\gamma|=1$, minimizing $\varE_\gamma$. 
We first claim that $\Om_\gamma\to B_{r_0}$ in $L^1(\R^d)$.  Indeed, clearly $e(0)\le \mathcal{E}_0(\Om_\gamma)\le \en(\Om_\gamma)$ for all $\gamma>0$, while on the other hand, 
$$ \limsup_{\gamma\to 0}\en(\Om_\gamma)\le \limsup_{\gamma\to 0}\en(B_{r_0}) = \mathcal{E}_0(B_{r_0}) = e(0).  $$
Thus, $\{\Om_\gamma\}_{\gamma>0}$ give minimizing sequences for $e(0)$.  
By Theorem~\ref{theorem:1} it follows that for every sequence $\gamma_n\to 0$ there is a subsequence along which $\{\Om_{\gamma_n}\}$ is compact in $L^1(\R^d)$, and converges to the minimizer $B_{r_0}$ of $e(0)$.  As the minimizer of the limit problem is unique, we conclude
 $\Om_\gamma\to B_{r_0}$ in $L^1(\R^d)$ as claimed.
      
Next, by Lemma \ref{lem:quasimin} $\Om_\gamma$ is an $\omega$-minimizer of the Euclidean perimeter in $O=\Rd\setminus \ol{B_\delta}$ for any $\delta>0$; hence, by Lemma \ref{lem:reg}(ii) $\pt^* \Om_\gamma \cap O$ is a $C^{1,\alpha}$-graph over the limit set $B_{r_0}$.  Note in particular that Lemma~\ref{lem:reg}(ii) assures that for any $\delta>0$, there exists $\gamma_\delta>0$ so that 
\beqn \label{annulus}
\pt^*\Om_\gamma \cap (B_{r_0-\delta}\setminus \overline{B_\delta}) = \emptyset
\eeqn
 for all $0<\gamma<\gamma_\delta$.  However,  since we cannot conclude quasi-minimality (hence, regularity) of $\Om_\gamma$ in the whole space we cannot {\it a priori} assume that $\pt\Om_\gamma \cap B_\delta$ is empty for all $\delta>0$.  Assume for a contradiction that for some $\delta>0$, 
 $\pt\Om_\gamma \cap B_\delta$ is non-empty; from \eqref{annulus} we may conclude that $|\Om_\gamma^c \cap B_\delta|>0$ for all $0<\gamma<\gamma_\delta$, and hence $\Om_\gamma$ is multiply-connected, and its boundary is disconnected into disjoint components. That is,
$\Om_\gamma = \tld{\Om}_\gamma \setminus \Om_\gamma^0$, where $\tld{\Om}_\gamma$ is simply connected, with smooth boundary, and
$\Om_\gamma^0 \defeq \Om_\gamma^c\cap B_\delta$, and thus 
$\pt^*\Omega_\gamma = \pt \tld{\Om}_\gamma \cup \pt^* \Omega_\gamma^0$, a disjoint union.
We define $m_\gamma \defeq |\Om_\gamma^0|$.  As $|\Omega_\gamma\triangle B_{r_0}|\to 0$ it follows that $m_\gamma\to 0$ and $|\tld{\Om}_\gamma| = 1 + m_\gamma$.   Moreover, let $R_\gamma = \left(\frac{1+m_\gamma}{\omega_d}\right)^{1/d}$ so that $|B_{R_\gamma}|=1+m_\gamma$.  In fact, using a Fuglede-type argument (cf. \cite{Fuglede}), we will prove that $m_\gamma=0$ for $\gamma$ sufficiently small, and therefore $\Omega_\gamma= B_{r_0}$ for all small $\gamma$.

Since $\pt \tld{\Om}_\gamma$ is a graph over $\pt B_{r_0}$ we can think of it as a graph over $\pt B_{R_\gamma}$ given by
	\[
		\pt \tld{\Om}_\gamma = \big\{  \tld{x} (1+\psi_\gamma(\tld{x})) \colon \tld{x}\in\pt B_{R_\gamma} \big\} 
	\]
for some $\psi_\gamma \in C^{1,\alpha}(\pt B_{R_\gamma})$ such that $\int_{\pt B_{R_\gamma}} \psi_\gamma \d\Hn =0$. Now, we can write both $\pt \tld{\Om}_\gamma$ and $\pt B_{R_\gamma}$ over the unit sphere $\Sd$. Namely,
	\[
		\pt \tld{\Om}_\gamma = \big\{ R_\gamma x (1+u_\gamma(x)) \colon x\in\Sd \big\} \quad \text{ and } \quad \pt B_{R_\gamma} =\big\{ R_\gamma x \colon x\in\Sd \big\},
	\]
where $u_\gamma(x) = \psi_\gamma(R_\gamma x)$. Note that $\int_{\Sd} u_\gamma \d \Hn = R_{\gamma}^{1-d} \int_{\pt B_{R_\gamma}} \psi_\gamma \d \Hn =0$. Computing the perimeter, we have
	\[
		\aPer{\tld{\Om}_\gamma} = R_\gamma^{d-1} \int_{\Sd} \a\big( R_\gamma x (1+u_\gamma(x))\big) (1+u_\gamma)^{d-1} \sqrt{1+\frac{|\nabla_{\tau} u_\gamma|^2}{(1+u_\gamma)^2}} \d \Hn,
	\]
where $\nabla_\tau$ denotes the gradient with respect to $\Sd$. Likewise,
	\[
		\aPer{B_{R_\gamma}} = R_\gamma^{d-1} \int_{\Sd} \a(R_\gamma x) \d \Hn.
	\]
Putting these two together,
\begin{align*}
	\aPer{\tld{\Om}_\gamma} - \aPer{B_{R_\gamma}} &= R_\gamma^{d-1} \int_{\Sd}  (1+u_\gamma)^{d-1} \a \big( R_\gamma x (1+u_\gamma(x))\big) \left[\sqrt{1+\frac{|\nabla_{\tau} u_\gamma|^2}{(1+u_\gamma)^2}}-1 \right]\d \Hn \\
												   &\qquad\qquad + R_\gamma^{d-1} \int_{\Sd} \Big[ (1+u_\gamma)^{d-1}\a\big(R_\gamma x(1+u_\gamma(x))\big) - \a(R_\gamma x) \Big] \d \Hn \\
												   &=:  R_\gamma^{d-1} I_1 + R_\gamma^{d-1} I_2.
\end{align*}

Now we will estimate $I_1$ and $I_2$. Let $\e>0$ be small so that $1+u_\gamma \geq 1/2$ and $\| u_\gamma \|_{C^1(\Sd)}\leq \e$ for $\gamma$ sufficiently small. Using the estimate $\sqrt{1+t} \geq 1+\frac{t}{2}-\frac{t^2}{8}$, we get that
\begin{align*}
	I_1 &\geq \int_{\Sd} R_\gamma^p(1+u_\gamma)^{p+d-1}\left[ \frac{1}{2}\frac{|\nabla_\tau u_\gamma|^2}{(1+u_\gamma)^2} -\frac{1}{8}\frac{|\nabla_\tau u_\gamma|^4}{(1+u_\gamma)^4} \right] \d \Hn \\
		 &\geq \left(\frac{1}{2} - C\e \right) \int_{\Sd} R_\gamma^p (1+u_\gamma)^{p+d-3} |\nabla_\tau u_\gamma|^2 \d \Hn \\
		 &\geq C \int_{\Sd} |\nabla_\tau u_\gamma|^2 \d \Hn
\end{align*} 
for $\e$ sufficiently small.

Since $\int_{\Sd} u_\gamma \d \Hn=0$ and $(1+t)^q \geq 1+qt$, we have
	\[
		I_2 = \int_{\Sd} R_\gamma^p \big( (1+u_\gamma)^{p+d-1}-1 \big) \d \Hn \geq 0.
	\]
These two estimates imply that
	\beqn \label{eq:per_diff}
		\aPer{\tld{\Om}_\gamma} - \aPer{B_{R_\gamma}} \geq C \| u_\gamma \|_{H^1(\Sd)}^2
	\eeqn
for some constant $C>0$ independent of $\gamma$.

As for the nonlocal term, \cite[Lemma 5.3]{FFMMM} (see also \cite[Remark 3.2]{FuPr}) implies that
	\[
		\N{B_{R_\gamma}} - \N{\tld{\Om}_\gamma} \leq C \| u_\gamma \|_{H^1(\Sd)}^2.
	\]
Furthermore,
	\beqn \label{eq:nonloc_diff}
		\V(B_{R_\gamma}\setminus \Om_\gamma^0)-\V(\Om_\gamma) \leq \V(B_{R_\gamma}) - \V(\tld{\Om}_\gamma) + Cm_\gamma \leq C \big( \| u_\gamma \|_{H^1(\Sd)}^2 + m_\gamma  \big).
	\eeqn
To see the first estimate above, let $\ol{R}>0$ be large enough such that $\Om_\gamma$, $B_{R_\gamma}\subset B_{\ol{R}}$ for all $\gamma>0$. For $\mathcal{I}(E,F)=\int_E\!\int_F |x-y|^{-\alpha}\d x \d y$, we have
	\[
		\V(B_{R_\gamma}\setminus \Om_\gamma^0)-\V(\Om_\gamma) = \V(B_{R_\gamma}) - \V(\tld{\Om}_\gamma) + 2 \mathcal{I}(\Om_\gamma^0,\tld{\Om}_\gamma)- 2 \mathcal{I}(\Om_\gamma^0,B_{R_\gamma}).
	\]
Since
	\begin{align*}
		\big|  \mathcal{I}(\Om_\gamma^0,\tld{\Om}_\gamma)- \mathcal{I}(\Om_\gamma^0,B_{R_\gamma})  \big| &\leq \int_{\Om_\gamma^0} \! \int_{\tld{\Om}_\gamma \triangle B_{R_\gamma}} \frac{1}{|x-y|^\alpha}\d x \d y \\
																																									&\leq |\Om_\gamma^0| \, \sup_{x\in B_{\ol{R}}} \int_{\tld{\Om}_\gamma \triangle B_{R_\gamma}} \frac{1}{|x-y|^\alpha}\d x \d y \leq C m_\gamma,
	\end{align*}
the first estimate in \eqref{eq:nonloc_diff} follows.

Again, using minimality of $\Om_\gamma$, i.e. $\varE_\gamma(\Om_\gamma) \leq \varE_\gamma(B_{R_\gamma} \setminus \Om_\gamma^0)$, and the lower and upper bounds \eqref{eq:per_diff} and \eqref{eq:nonloc_diff} we get that
	\begin{align*}
		C \| u_\gamma \|_{H^1(\Sd)}^2 &\leq \aPer{\tld{\Om}_\gamma} - \aPer{B_{R_\gamma}} = \aPer{\Om_\gamma} - \aPer{B_{R_\gamma}\setminus\Om_\gamma^0} \\
													  &\leq \gamma\, \big( \V(B_{R_\gamma}\setminus \Om_\gamma^0) - \V(\Om_\gamma) \big) \leq C\gamma\, \big( \| u_\gamma \|_{H^1(\Sd)}^2 + m_\gamma \big).
	\end{align*}
Thus, for $\gamma$ sufficiently small,
	\beqn \label{eq:H1_est}
		\| u_\gamma \|_{H^1(\Sd)}^2 \leq C\gamma \, m_\gamma.
	\eeqn
This, in turn, implies that
	\beqn \label{eq:en_upper_bd}
		\varE_\gamma(B_{R_\gamma}\setminus \Om_\gamma^0)-\varE_\gamma(\Om_\gamma) \leq C\gamma\,m_\gamma.
	\eeqn

The last step in the argument is a lower bound on the above difference.  In fact, we will show that replacing the domain $\Om_\gamma=\tld{\Omega}_\gamma\setminus \Omega_\gamma^0$ by the sphere $B_{r_0}$ results in a much larger energy difference, due to the reduction of the radius $R_\gamma$ to $r_0$.  We observe this right away, in estimating the difference in perimeter, $\aPer{B_{R_\gamma}}+\aPer{\Om_\gamma^0}-\aPer{B_{r_0}}$. Let $\rho_\gamma = (m_\gamma/\omega_d)^{1/d}$ and note that
	\[
		R_\gamma =(r_0^d+\rho_\gamma^d)^{1/d} = r_0 \left( 1+\left(\frac{\rho_\gamma}{r_0}\right)^d \right)^{1/d}.
	\]
Now, for $m_\gamma$ (hence, for $\rho_\gamma/\omega_d$) small, we have
	\begin{align*}
		\aPer{B_{R_\gamma}}+\aPer{\Om_\gamma^0}-\aPer{B_{r_0}} &\geq \aPer{B_{R_\gamma}}+\aPer{B_{\rho_\gamma}}-\aPer{B_{r_0}} \\
																								    &= d\omega_d (R_\gamma^{d-1+p}+\rho_\gamma^{d-1+p}-r_0^{d-1+p}) \\
																								    &= d\omega_d r_0^{d-1+p} \left( \left[ 1+\left(\frac{\rho_\gamma}{r_0}\right)^d \right]^{(d-1+p)/d} + \left(\frac{\rho_\gamma}{r_0}\right)^{d-1+p} -1 \right) \\
																								    &=d\omega_d r_0^{p-1} \left( 1 + \frac{p-1}{d} \right) \rho_\gamma^d + O(\rho_\gamma^{d-1+p})\\
																								    &= \bar{C}\,m_\gamma + o(m_\gamma).
	\end{align*}
Combining this estimate with \eqref{eq:nonloc_diff} and \eqref{eq:H1_est}, and recalling that the ball maximizes the Riesz potential, for sufficiently small $\gamma$, we get
	\beqn \label{eq:en_lower_bd}
		\begin{aligned}
			\varE_\gamma(B_{R_\gamma}\setminus \Om_\gamma^0)-\varE_\gamma(B_{r_0}) &= \aPer{B_{R_\gamma}}+\aPer{\Om_\gamma^0}-\aPer{B_{r_0}} + \gamma \big( \V(B_{R_\gamma}\setminus\Om_\gamma^0)-\V(B_{r_0}) \big) \\
																																			   &\geq (\bar{C}-\gamma C) m_\gamma \geq C m_\gamma.
		\end{aligned}
	\eeqn
Now \eqref{eq:en_upper_bd} and \eqref{eq:en_lower_bd} imply that
	\[
		Cm_\gamma + \varE_\gamma(B_{r_0}) \leq \varE_\gamma(B_{R_\gamma}\setminus\Om_\gamma^0) \leq \varE_\gamma(\Om_\gamma) + C\gamma \, m_\gamma \leq \varE_\gamma(B_{r_0}) + C\gamma\, m_\gamma.
	\]	
Hence, for $\gamma$ sufficiently small, $Cm_\gamma \leq 0$, i.e., $m_\gamma=0$. This implies that $\Om_\gamma^0=\emptyset$ for $\gamma$ small and that $\pt \Om_\gamma$ is a graph over $\pt B_{r_0}$.

Therefore, running through the Fuglede-type argument one more time, we get
	\[
		C \| u_\gamma \|_{H^1(\Sd)}^2 \leq \aPer{\Om_\gamma} - \aPer{B_{r_0}} \leq \gamma\big(\N{B_{r_0}} - \N{\Om_\gamma}\big)\leq C\gamma\,\| u_\gamma \|_{H^1(\Sd)}^2.
	\]
Hence, for $\gamma$ sufficiently small $u_\gamma \equiv 0$, i.e., $\Om_\gamma = B_{r_0}$.
\end{proof}

%%%%%%%%%%%%%%%%%%%%%%%%%%%%%%%%%%%%%%%%%%%%%%%%%%%%%%%%%%%%%%%%%
%%%%%%%%%%%%%%%%%%%%%%%%%%%%%%%%%%%%%%%%%%%%%%%%%%%%%%%%%%%%%%%%%

\subsection*{Acknowledgments}
The authors would like to thank Marco Bonacini and Gian Paolo Leonardi for valuable discussions regarding the properties of density perimeters. We also thank the reviewer for carefully reading the paper and providing many useful suggestions.  SA, LB, and AZ were supported via an NSERC (Canada) Discovery Grant. AZ was also partially funded by ANID Chile under grants Becas Chile de Postdoctorado en el Extranjero N$^{\circ}$ 74200091 and FONDECYT de Iniciaci\'on en Investigaci\'on N$^{\circ}$ 11201259.

\bibliographystyle{IEEEtranS}
\def\url#1{}
\bibliography{references}

\end{document}